\documentclass[journal]{IEEEtran}
\usepackage{color}

\usepackage{amssymb}
\usepackage{comment}

\usepackage{bbm}
\usepackage{mathrsfs}
\usepackage{verbatim}
\usepackage{amsmath}
\usepackage{amsfonts}
\usepackage{bm}
\usepackage{graphicx}
\usepackage{caption}
\usepackage{float}

\usepackage{amsthm}
\usepackage{mathtools}
\usepackage{cite}
\usepackage{bigints}
\usepackage{url}
\usepackage{hyperref}

\hypersetup{
    colorlinks=true,
    linkcolor=blue,
    citecolor =black,
    filecolor=black,      
    urlcolor=blue,
    }
    
\usepackage{subfig}
\usepackage{savesym}

\usepackage{algorithm}
\usepackage{algorithmicx} % <===========================================
\usepackage{algpseudocode} %
\savesymbol{AND}

\usepackage{bigints}

\usepackage[group-separator={,},group-minimum-digits={3}]{siunitx}

\theoremstyle{definition}
\newtheorem{assumption}{Assumption}
\newtheorem{definition}{Definition}
\newtheorem{remark}{Remark}
\newtheorem{problem}{Problem}

\theoremstyle{plain}

\newtheorem{theorem}{Theorem}
\newtheorem{lemma}{Lemma}
\newtheorem{corollary}{Corollary}

\usepackage{csquotes}
\MakeOuterQuote{"}
\usepackage{booktabs}
\usepackage{array}   % for \newcolumntype macro

\makeatletter

\newcommand{\Rmnum}[1]{\expandafter\@slowromancap\romannumeral #1@}
\allowdisplaybreaks[0]

\begin{document}
       
\title{Optimal Control of Connected and Automated Vehicles at Multiple Adjacent Intersections}

\author{Behdad Chalaki, {\itshape{Student Member, IEEE}} and Andreas A. Malikopoulos, {\itshape{Senior Member, IEEE}}
\thanks{This research was supported in part by ARPAE's NEXTCAR program under the award number DE-AR0000796 and by the Delaware Energy Institute (DEI).}%
\thanks{The authors are with the Department of Mechanical Engineering, University of Delaware, Newark, DE 19716 USA (email: \texttt{bchalaki@udel.edu; andreas@udel.edu).}} }

%\markboth{IEEE Transactions on Control Systems Technology}{Chalaki and Malikopoulos : Optimal Control of Connected and Automated Vehicles at Multiple Adjacent Intersections}

\maketitle

\begin{abstract}     
In this paper, we establish a decentralized optimal control framework for connected and automated vehicles (CAVs) crossing multiple adjacent, multi-lane signal-free intersections to minimize energy consumption and improve traffic throughput. Our framework consists of two layers of planning. In the upper-level planning, each CAV computes its optimal arrival time at each intersection recursively along with the optimal lane to improve the traffic throughput. In the low-level planning, we formulate an energy-optimal control problem with interior-point constraints, the solution of which yields the optimal control input (acceleration/deceleration) of each CAV to cross the intersections at the time specified by the upper-level planning. Moreover, we extend the results of the proposed bi-level framework to include a bounded steady-state error in tracking the optimal position of the CAVs. Finally, we demonstrate the effectiveness of the proposed framework through simulation for symmetric and asymmetric intersections and comparison with traditional signalized intersections.
\end{abstract}

\section{Introduction}

\IEEEPARstart{O}{ver} the last few decades, the urban population of the world has proliferated. Today, $55$\% of the world's population lives in urban areas; this ratio is expected to increase to $68$\% by $2050$ \cite{united20182018}. However, in urban areas, road capacity has not grown at the same pace, resulting in traffic congestion. Traffic congestion has been persistently growing from $1982$ to $2017$ in US urban areas \cite{Schrank2019}. 
In addition, traffic congestion has a negative impact on traffic safety. In $2018$, there were 6M traffic accidents in the US resulting in more than $35$K fatalities and $2.5$M people injured \cite{USDOT2}.

One of the promising ways to mitigate traffic congestion and improve safety is integrating information and communication technologies in cities by utilizing connected and automated vehicles (CAVs) \cite{Klein2016a,Melo2017a}. After the seminal work of Levine and Athans \cite{Levine1966,Athans1969} on safely coordinating vehicles at merging-road ways, several research efforts have explored the benefits of coordinating CAVs in traffic scenarios, such as urban intersections, merging roadways, and speed reduction zones to eliminate congestion in a transportation network while preserving safety.

\subsection{Related Work}
Over the last few years, there has been an increased interest in investigating different approaches for coordination of CAVs at intersections. These approaches can be categorized into two major groups, namely, centralized and decentralized. In centralized approaches, there is at least one task in the system that is globally decided for all CAVs by a single central controller. In decentralized approaches, CAVs are treated as autonomous agents that optimize specific performance criteria (e.g., fuel efficiency, travel time) through vehicle to vehicle (V2V) and/or vehicle to infrastructure (V2I) communication. 

Several studies have developed a centralized approach to improve travel time while guaranteeing safety at intersections \cite{Dresner2008,Lee2012a,gregoire2014priority,fayazi2018mixed}. Dresner and Stone \cite{Dresner2008} introduced a reservation-based scheme which requires CAVs to reserve a space-time slot inside the intersection. Lee and Park \cite{Lee2012a} minimized the total length of overlapped trajectories of CAVs crossing an intersection. Gregoire \textit{et al.} \cite{gregoire2014priority} decomposed the coordination problem into a central priority assignment and trajectory planning. Given the priority assignment, they planned a safe trajectory with either maximum or minimum control inputs. Fayazi and Vahidi \cite{fayazi2018mixed} proposed a framework based on the arrival time of CAVs at an intersection. Then, they converted the arrival time scheduling problem to a central mixed-integer-linear-program (MILP). 

Other efforts have also considered improving fuel efficiency in the coordination problem. Bichio and Rakha \cite{bichiou2018developing} demonstrated an improvement in fuel efficiency of CAVs by jointly minimizing the travel time and control efforts for $M\in\mathbb{N}$ closest CAVs to the intersection; however, the proposed method is not real-time implementable (e.g., $2$-$5$ minutes for $M=4$). Borek \textit{et al.} \cite{borek2019economic} presented the energy-optimal control for heavy-duty trucks, which employs an online model predictive control (MPC) to track the optimal solution obtained off-line using dynamic programming. Du \textit{et al.} \cite{du2018hierarchical} presented a three-layered hierarchical coordination strategy for CAVs at multiple intersections. In the top layer, each controller generates the desired road speed to balance the traffic density over multiple intersections. In the next layer, each controller computes a trajectory for each CAV minimizing the deviation from the desired road speed and satisfying lateral safety in intersections. In the last layer, each CAV uses MPC to track the prescribed trajectory while avoiding rear-end collision. 

There has also been a series of papers to investigate overriding CAVs at intersections only if their current inputs lead to a collision \cite{Colombo2015,Colombo2014,Ahn2014,Ahn2016a,DeCampos2015a}. Colombo and Del Vecchio \cite{Colombo2015} demonstrated that checking whether current inputs of CAVs lead to a collision is equivalent to a scheduling problem, the solution of which yields constructing a least restrictive controller to ensure maintaining the state of the system within the maximal controlled invariant set. Colombo \cite{Colombo2014} extended the results to a network with arbitrarily many intersections by decoupling intersections and handling them in isolation.

To date, several research efforts in the literature have explored decentralized approaches for the coordination of CAVs at intersections. One of the early efforts was proposed by Makarem and Gillet \cite{Makarem2012} using a navigation function aimed at minimizing the energy consumption of each CAV. Wu \textit{et al.} \cite{wu2014distributed} presented an algorithm based on a mutual exclusion in which CAVs compete for the privilege of crossing the intersection through V2V communication. Focusing on V2V communication, Azimi \textit{et al.} \cite{azimi2014stip} proposed various intersection protocols to improve traffic throughput while avoiding collisions. In their approach, the control zone is considered a grid divided into small cells. In their most advanced protocol, when there is a potential conflict at a cell, a lower-priority CAV can either cross the conflicting cell or arrive at the cell after the higher-priority CAV has exited the cell. 
Hult \textit{et al.} \cite{hult2018optimal} presented a bi-level coordination scheme for CAVs crossing a signal-free intersection. The first level is a central controller, which yields the required arrival times for CAVs at the intersection along with CAVs' crossing orders. In the second level, the authors considered a local MPC given the arrival time computed from the first level. Other authors have also formulated a local MPC problem for each CAV with defined crossing priorities \cite{kim2014mpc,campos2014cooperative,kloock2019distributed}.

In decentralized coordination of CAVs, various research efforts have been reported in the literature using optimal control techniques to find closed-form solutions. Malikopoulos \textit{et al.} \cite{Malikopoulos2017} established a decentralized coordination framework for CAVs at an intersection which consists of throughput maximization and energy minimization problems. In the throughput maximization problem, each CAV computes its arrival time at the merging zone, i.e., the area of potential lateral collisions, based on a first-in-first-out (FIFO) queuing policy. By restricting CAVs to have constant speed at the merging zone, each CAV derives its energy-optimal control input from the entry of the control zone until it reaches the merging zone considering speed and control constraints. In sequel papers, Malikopoulos and Zhao \cite{malikopoulos2019ACC} enhanced the decentralized framework for a single intersection by including speed-dependent rear-end safety constraint, while Mahbub \textit{et al.} \cite{Mahbub2019ACC} presented the unconstrained solution for two adjacent intersections.
In a different approach for coordination of CAVs at an intersection, the objective function of each CAV was formulated by jointly minimizing travel time and energy consumption, where the analytical solution was presented in \cite{zhang2019decentralized} with minimum distance rear-end safety constraint, and in \cite{zhang2019joint} with speed-dependent rear-end safety constraint. 

 A comprehensive discussion of the research efforts that has been reported in the literature to date in control and coordination of CAVs is provided in \cite{guanetti2018control} and \cite{Rios-Torres2017}.

\subsection{Contributions of This Paper}

Although there have been many studies reporting on the coordination of CAVs at signal-free intersections, coordination of CAVs at multiple adjacent intersections has been assessed to a very limited extent. 
One of the main drawbacks of considering each intersection in isolation is neglecting the effects of the downstream intersection on the upstream intersection. Mahbub \textit{et al.} \cite{mahbub2020decentralized} presented coordination of CAVs at a traffic corridor consisting of multiple traffic scenarios by considering each scenario in isolation. Extending the single intersection results to the multiple intersections might be inefficient and sub-optimal. For example, applying FIFO queuing policy \cite{Malikopoulos2017,zhang2019decentralized} to find the sequence of CAVs to enter the merging zone in a single intersection will result in unnecessary slowdowns of the CAVs at multiple intersections \cite{chalaki2020TITS}.

In an earlier work \cite{chalaki2020TITS}, we presented a hierarchical control framework for two adjacent single-lane intersections. In the upper-level, we formulated a scheduling problem for each CAV to compute the arrival times at each zone, the solution of which minimized the total travel time of the CAV, and it was solved using a MILP. In the low-level problem, for each zone, the decentralized energy optimal control problem was formulated, and closed-form solutions were provided. However, in \cite{chalaki2020TITS}, the focus was on minimizing the travel-time rather than energy consumption, and thus the energy consumption was not improved significantly compared to the baseline scenario with traffic signals. In addition, since the approach proposed in \cite{chalaki2020TITS} employs a MILP in the upper-level layer, it may not be appropriate to investigate multiple intersections, due to computational complexity.

In this paper, we present a bi-level decentralized coordination framework for CAVs at multiple adjacent, multi-lane signal-free intersections that are closely distanced. In the upper-level planning, each CAV recursively computes the energy-optimal arrival time at each intersection along its path, while ensuring both lateral and rear-end safety. In the low-level planning, we formulate an optimal control problem for each CAV with interior-point constraints, the solution of which yields the energy optimal control input, given the time from the upper-level problem. The contributions of this paper are: (1) the development of a bi-level optimization framework to coordinate CAVs at multiple adjacent, multi-lane intersections aimed at decreasing stop-and-go driving and fuel consumption; (2) the enhancement of the upper-level planning layer to include the lane-changing maneuver to improve the traffic throughput; (3) a complete analysis of the low-level optimization problem including interior-point constraints; and (4) the enhancement of the bi-level framework to account for a bounded steady-state error in tracking the positions of CAVs.

 \subsection{Comparison with Related Work}

The proposed framework advances the state of the art in the following ways. First, rather than investigating intersections in isolation \cite{Colombo2014,Rios-Torres2017,mahbub2020decentralized,Malikopoulos2017}, we propose a decentralized framework by considering the effects of intersections' interdependence methodically.
Second, we ensure lateral safety through a decentralized upper-level planning, as opposed to a strict FIFO queuing policy \cite{Malikopoulos2017,zhang2019decentralized,zhang2019joint,bichiou2018developing}, or a centralized controller  \cite{Dresner2008,Lee2012a,xu2019cooperative,kamal2014vehicle,du2018hierarchical,guney2020scheduling,Gregoire2014a,hult2018optimal}. Third, in contrast to the research efforts reported in the literature to date, in the upper-level planning, we allow lane-changing maneuvers. Finally, our bi-level framework is enhanced to guarantee safety in the presence of a bounded steady-state error in tracking the positions of CAVs. 

\subsection{Organization of This Paper}
The rest of the paper is organized as follows. 
In Section \ref{section2}, we introduce the modeling framework, while in Section \ref{section3} and \ref{section4}, we provide the upper-level planning and low-level planning with their corresponding solutions, respectively. In Section \ref{section5}, we enhance the framework to include the deviation from the nominal planned position. Finally, we provide simulation results in Section \ref{section6} and offer concluding remarks along with a discussion for a future research direction in Section \ref{section7}.

\section{Problem Formulation}\label{section2}
\begin{figure}[ht]
\centering
\includegraphics[width=0.95\linewidth]{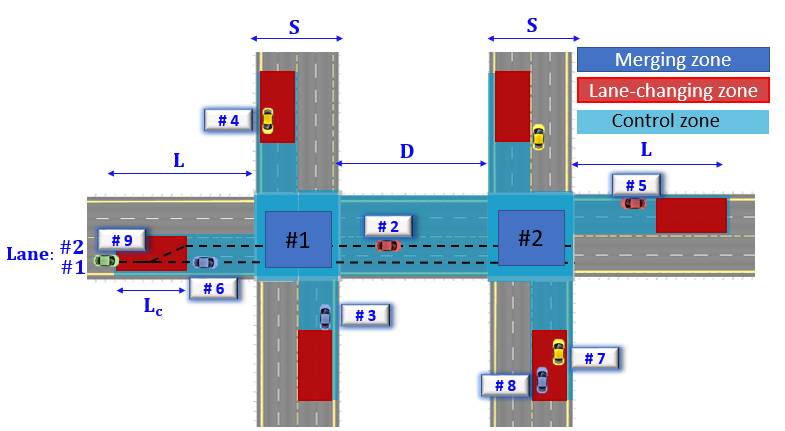}
\caption{Bird-eye view of two interconnected intersections.}
\label{fig:mlLanechange4}%
\end{figure}

We consider multiple adjacent intersections closely distanced from each other (see Fig. \ref{fig:mlLanechange4} for two adjacent intersections). There is a \textit{coordinator} that stores information about all intersections' geometric parameters and the planned trajectories of CAVs. The coordinator does not make any decision and it only acts as a \textit{\mbox{{database}}} among the CAVs. We define  a \textit{{contol zone}} in which 
the coordinator can communicate with the CAVs traveling inside the control zone. We call the areas inside the control zone where lateral collisions may occur \textit{{merging zones}}.

\begin{definition}
The set of merging zones indexed uniquely in the control zone, is given by $\mathcal{Z} \coloneqq\{1,\dots,n_z\},~n_z\in\mathbb{N}$, where $n_z$ is the total number of merging zones in the control zone.
\end{definition}
\begin{definition}
The set of all same-directional lanes across all roads that are connected to the intersection, is given by $\mathcal{L}\coloneqq\{1,\dots,n_l\},~ n_l\in\mathbb{N}$, where $1$ and $n_l$ are indices for the rightmost lane and leftmost lane, respectively (Fig. \ref{fig:mlLanechange4}).
\end{definition}

Without loss of generality, in our analysis, we consider multiple intersections with two lanes on each road that is connected to an intersection, i.e., $\mathcal{L}\coloneqq\{1,2\}$. We consider the distance between the two adjacent intersections to be $D$, and a merging zone of length $S = 4w$, where $w\in\mathbb{R}^{+}$ is the lane width. In this paper, we limit our analysis to the cases that no left/right turns are allowed, and lane-change maneuver is only possible inside the lane-changing zone which we describe next.

%\subsection{Local Coordination System}

Let $N(t)\in\mathbb{N}$ be the total number of CAVs that have entered the control zone by the time $t\in\mathbb{R}^{+}$, and $\mathcal{N}(t)=\{1,\ldots,N(t)\}$ be the queue that designates the order that each CAV entered the control zone. Upon entering the control zone, the coordinator assigns each CAV an integer index equal to $N(t)+1$. If two or more CAVs enter the control zone simultaneously, the CAV with the shorter path is assigned lower position in the queue; however, if the length of their paths is the same, then their positions are randomly chosen by the coordinator. 
The coordinator removes any CAV from $\mathcal{N}(t)$ when they exit the control zone. When there is no CAV inside the control zone, then $\mathcal{N}(t)=\emptyset$.

We model the dynamics of each CAV $i\in\mathcal{N}(t)$ as a double integrator, i.e.,
\begin{align}
\begin{aligned}\label{27a}
\dot{p}_i(t)=v_i(t),\\
\dot{v}_i(t)=u_i(t),
\end{aligned}\end{align}
where $p_{i}(t)\in\mathcal{P}_{i}$, $v_{i}(t)\in\mathcal{V}_{i}$, and
$u_{i}(t)\in\mathcal{U}_{i}$ denote position, speed, and acceleration at $t\in\mathbb{R}^{+}$, respectively. Let $\mathbf{x}_{i}(t)=\left[p_{i}(t), v_{i}(t)\right]^\top$ and $u_{i}(t)$ be the state
and control input of the CAV $i$ at time $t$, respectively, $t_{i}^{0}\in\mathbb{R}^{+}$ be the time that CAV $i\in\mathcal{N}(t)$
enters the control zone, and $t_{i}^{f}>t_i^0\in\mathbb{R}^{+}$ be the time that CAV $i$ exits the control zone (the merging zone of the last intersection along its path).
For each CAV $i\in\mathcal{N}(t)$ the control input and speed are bounded by 
\begin{equation}\label{uconstraint}
    u_{i,\min}\leq u_i(t)\leq u_{i,\max},
\end{equation}
\begin{equation}\label{vconstraint}
    0\leq v_{\min}\leq v_i(t)\leq v_{\max},
\end{equation}
where $u_{i,\min},u_{i,\max}$ are the minimum and maximum control inputs and $v_{\min},v_{\max}$ are the minimum and maximum speed limit, respectively. The sets $\mathcal{P}_{i}$,
$\mathcal{V}_{i}$, and $\mathcal{U}_{i}$, $i\in\mathcal{N}(t),$
are complete and totally bounded subsets of $\mathbb{R}$.

\begin{definition} The lane-changing zone $\Lambda$ is the interval with length $L_c$ located at the entry of the control zone, where CAVs can change lanes (Fig. \ref{fig:mlLanechange4}), i.e.,
 \begin{equation}
     \Lambda \coloneqq [p_i(t_i^0),~p_i(t_i^0)+L_c]\subset \mathcal{P}_{i} ,~ i\in\mathcal{N}(t). 
 \end{equation}
\end{definition}
\begin{definition}
The lane-changing occupancy interval $\Gamma_i$ is the time interval that CAV $i\in\mathcal{N}(t)$ occupies the lane-changing zone, i.e.,
\begin{equation}
    \Gamma_i \coloneqq \left\{t~|~ t>t_i^0,~ p_i(t) \leq L_c,~t\in\mathbb{R}^+\right\}.
\end{equation}
\end{definition}

\begin{definition}
For each CAV $i\in\mathcal{N}(t)$, $l_i^0,~l_i^f\in\mathcal{L}$ denote the lane that CAV $i$ occupies before and after the lane-changing zone, respectively.
\end{definition}

\begin{definition}\label{Defn:Z_i}
 Let $\mathcal{Z}_i \coloneqq \{z_1,\dots,z_n\}\subseteq\mathcal{Z}$ be the set of merging zones that CAV $i\in\mathcal{N}(t)$ crosses, while $z_1$ and $z_n\in\mathcal{Z}_i,~n\in\mathbb{N}$, denote the first and $n$'th merging zone that CAV $i$ crosses along its path.
\end{definition}

\begin{definition}
For each CAV $i\in\mathcal{N}(t)$, we denote $t_i^{z^\ast}\in\mathbb{R}^{+}$ and $l_i^{f^\ast}\in\mathcal{L}$ to be the optimal arrival time at the entry of zone $z\in\mathcal{Z}_i$ and the optimal lane occupied after the lane-changing zone, respectively.
\end{definition}

For CAV $i\in\mathcal{N}(t_i^0)$, CAV $j\in\mathcal{N}(t_i^0)$, $j<i$, belongs to one of the following time-invariant subsets, determined upon CAV $i$'s entrance at the control zone:
\begin{enumerate}
    \item $\mathcal{A}_i^l, l\in\mathcal{L}$, is the set of all CAVs that travel on lane $l$ after a lane-changing zone with the same direction and destination as CAV $i$.
    \item $\mathcal{B}_i^{z}, z\in\mathcal{Z}_i$, is the set of all CAVs which may cause collision with CAV $i$ at the merging zone $z$ at time $t \geq t_i^0$. 
    \item $\mathcal{C}_i$, is the set of all CAVs with a different origin-destination pair from CAV $i$, without any potential conflict for $t\geq t_i^0$.
\end{enumerate}
For example, consider CAV \#$9$ in Fig. \ref{fig:mlLanechange4}. We have: $\mathcal{A}_9^1=\{6\}$, $\mathcal{A}_9^2=\{2\}$, $\mathcal{B}_9^{1}=\{3,4\}$, $\mathcal{B}_9^{2}=\{7,8\}$, $\mathcal{C}_9=\{5\}$.

To ensure the absence of rear-end collision between CAV $i\in\mathcal{N}(t)$ and a preceding CAV $k\in\mathcal{N}(t)$, we impose the following rear-end safety constraint
\begin{equation}\label{RearEndCons}
 p_k(t)-p_i(t)\geq \delta,
\end{equation}
where $\delta\in\mathbb{R}^+$ is a constant safe distance. Note that, since we study urban intersections, the average speed variation is not significant, thus considering a constant safe distance is reasonable. However, for scenarios in which the average speed variation is not negligible, one may consider speed-dependent safe distance discussed in \cite{malikopoulos2019ACC,chalaki2020TITS,Malikopoulos2020,zhang2019joint}. 

\begin{definition}
The lane changing maneuver for CAV $i\in\mathcal{N}(t)$ is defined to be feasible, if upon arriving at the lane-changing zone, no other CAV occupies the lane-changing zone, i.e., $\{t_i^0\} \cap\Gamma_j = \emptyset$ for all $j\in\mathcal{A}_i^l$, $l\in\mathcal{L}$. 
\end{definition}

\begin{definition}
For each CAV $i\in\mathcal{N}(t)$, $\mathcal{T}^z_i$ is the set of optimal arrival times of CAVs at zone $z\in\mathcal{Z}_i$ that belong to $\mathcal{B}_i^z$, i.e.,
\begin{equation}
    \mathcal{T}^z_i \coloneqq \{t_j^{z^\ast}~|~ j\in\mathcal{B}_i^{z}\}.
\end{equation}
\end{definition}

In the modeling framework described above, we impose the following assumptions:

\begin{assumption}\label{feassibleAssum}
Upon entering the control zone, CAV $i\in\mathcal{N}(t)$ has at least a distance $\delta\in\mathbb{R}^+$ from any preceding CAV traveling at any lane $l\in\mathcal{L}$.
\end{assumption}
\begin{assumption}\label{feassibleAssum2}
None of the speed constraints is active for each CAV $i\in \mathcal{N}(t)$ at the entry of the control zone.
\end{assumption}
Assumptions \ref{feassibleAssum} and \ref{feassibleAssum2} are imposed to ensure that the initial state is feasible. These are reasonable assumptions since CAVs are automated, and so there is no compelling reason for them to activate any of the state constraints by the time they enter the control zone.

\section{Upper-level Planning}\label{section3}
The objective of each CAV at the entry of the control zone is to derive the optimal control input (acceleration/deceleration) aimed at minimizing fuel consumption and improving traffic throughput by eliminating stop-and-go-driving. To achieve this aim, we establish a decentralized control framework consisting of two layers of planning. In the upper-level planning, each CAV $i\in\mathcal{N}(t)$ recursively computes the arrival time at each merging zone along its path with the optimal lane to occupy after lane-changing zone in order to improve traffic throughput and energy consumption. The outputs of the upper-level planning become the inputs for the low-level planning, which we describe in Section \ref{section4}.
\begin{definition}
For each CAV $i\in\mathcal{N}(t)$, $v_{\text{avg}}^z$ is its average speed inside the merging zone $z\in\mathcal{Z}_i$, i.e.,
\begin{equation}
   v_{\text{avg}}^z = \dfrac{\int_{t_i^{z}}^{t_i^{z}+\Delta t_i^z} v_{i}(t)dt}{\Delta t_i^z},
\end{equation}
where $t_i^{z}$, $\Delta t_i^z\in\mathbb{R}^+$ are the arrival time at zone $z$ and the time that it takes for CAV $i$ to exit the merging zone $z$, respectively.
\end{definition}

To increase safety and improve the throughput of CAV $i\in\mathcal{N}(t)$ while traveling inside the merging zone $z\in\mathcal{Z}_i$, we impose a constant average speed inside the merging zone equal to the speed that CAV $i$ entered the control zone ( $v_{\text{avg}}^z\coloneqq v_i(t_i^0)$). This results in traveling at the merging zone with constant time $\Delta t_i^z = \dfrac{\Delta x_i^z}{v_{\text{avg}}^z}$, where $\Delta x_i^z$ is the distance traveled at merging zone $z$ for CAV $i$. As mentioned earlier, since no turn is allowed inside the merging zone, $\Delta x_i^z$ is the same for all CAVs, and it is equal to $S=4w$ (see Fig. \ref{fig:mlLanechange4}). It should be noted that, imposing desired average speed inside the merging zone is different from setting a constant speed as in \cite{Malikopoulos2017}, and thus it is less restrictive since CAV's speed can vary inside the merging zone as long as it satisfies the desired average speed. Moreover, to minimize the energy consumption of CAV $i\in\mathcal{N}(t)$ inside the control zone, we minimize transient engine operation, $L^2$-norm of the control input in $[t_i^0,~t_i^f]$, which was shown to have direct benefit in fuel consumption and emission \cite{Malikopoulos:2013aa}.

\begin{lemma}\label{lem:tbar}
The arrival time of CAV $i\in\mathcal{N}$ at the merging zone of $n$'th intersection $z_n\in\mathcal{Z}_i$ along its path, without considering rear-end or lateral safety constraints, minimizing the energy-consumption is denoted by $\bar{t}_i^{z_n}$ and is computed recursively as follows
\begin{equation}\label{tbar}
    \bar{t}_i^{z_n} = \begin{cases}
            t_i^0+\dfrac{L}{v_i(t_i^0)}  ,&\text{for}\;\: z_1,\\
            \\
            t_i^{z_{n-1}^\ast}+\Delta t_i^{z_{n-1}}+\dfrac{D}{v_i(t_i^0)},&\text{otherwise.}
          \end{cases}
\end{equation} 
\end{lemma}

\begin{proof}
There are two cases to consider: Case 1: $z_1$ and Case 2: $z_n$.

Case 1: Suppose CAV $i\in\mathcal{N}(t)$ enters the control zone at $t=t_i^0$, and let the arrival time at merging zone of the first intersection $z_1\in\mathcal{Z}_i$ along its path be $t_i^{z_1}$. 
Let $\bar{t}_i^{z_1}$ be the arrival time at $z_1\in\mathcal{Z}_i$ minimizing the following cost function $J_i(u_i(t),t_i^{z_1})= {\dfrac{1}{2}} \bigints_{t_i^{0}}^{t_i^{z_1}} u_{i}(t)^2dt$, without considering rear-end safety or lateral safety constraints.
For the unconstrained case, the Hamiltonian is
\begin{equation}
    H_i(t,p_i(t),v_i(t),u_i(t))=\frac{1}{2}u_i(t)^2+\lambda_i^p v_i (t)+\lambda_i^v u_i(t),
\end{equation}
where \(\lambda_i^p\) and \(\lambda_i^v\) are costates. Applying the Euler-Lagrange optimality conditions, the optimal control input minimizing the cost function $J_i(u_i(t),t_i^{z_1})$ is $u^\ast_i(t) = -{\lambda_i^v}^\ast= a_it+b_i$ \cite{Malikopoulos2017}, where $a_i$ and $b_i$ are constants of integration. We also have ${\lambda_i^p}^\ast= a_i$. Since the speed at $t= \bar{t}_i^{z_1}$ is not specified, we have $\lambda_i^v(\bar{t}_i^{z_1}) = 0$. In addition, since $t_i^{z_1}$ is not defined, we have the following transversality condition
$H_i(\bar{t}_i^{z_1},p^\ast_i(t),v^\ast_i(t),u^\ast_i(t))=0$. From the transversality condition, we have 
$H_i(\bar{t}_i^{z_1},p^\ast_i(t),v^\ast_i(t),u^\ast_i(t))= {\lambda_i^p}^\ast v^\ast_i(\bar{t}_i^{z_1})=0$, and because $v^\ast_i(\bar{t}_i^{z_1})\neq0$, we get 
\begin{equation}\label{eq:zeroAcc}
  {\lambda_i^p}^\ast = 0 \Rightarrow u^\ast_i(t)=0 ,~\forall~ t \in[t^0_i,\bar{t}_i^{z_1}].
\end{equation}
Hence, CAV $i$ cruises with $v_i(t_i^0)$, and $\bar{t}_i^{z_1} = t_i^0+\dfrac{L}{v_i(t_i^0)}$.
 
Case 2: For the $n$'th merging zone, $z_n\in\mathcal{Z}_i$, let us consider the optimal arrival time at the upstream merging zone, $z_{n-1}\in\mathcal{Z}_i$, to be $t_i^{z_{n-1}^\ast}$. Let $\bar{t}_i^{z_n}$ be the energy-efficient arrival time at $z_n\in\mathcal{Z}_i$ without considering rear-end safety or lateral safety constraint. As it was shown in Case $1$, by neglecting the safety constraint, CAV $i\in\mathcal{N}(t)$ cruises with $v_i(t_i^0)$ to minimize the energy consumption. Thus, $\bar{t}_i^{z_n} = t_i^{z_{n-1}^\ast}+\Delta t_i^{z_{n-1}}+\dfrac{D}{v_i(t_i^0)} $. \end{proof} 

\begin{remark}
To consider the impact of the upstream merging zone $z_{n-1}$ on the merging zone $z_n$ for CAV $i\in\mathcal{N}(t)$, we need a recursive formulation to relate the arrival time at $z_n$ to the optimal arrival time at $z_{n-1}$.
\end{remark}

In order for CAV $i\in\mathcal{N}(t)$ to avoid the lateral collision with CAV $j\in\mathcal{B}_i^{z}$, it can either arrive at merging zone $z\in\mathcal{Z}_i$ after CAV $j$ exits the merging zone $z$, or exit the merging zone $z$ before CAV $j$ enters the merging zone $z$. This is formulated as either
\begin{equation} \label{eq:B_iz}
t_i^{z^\ast} \geq t_j^{z^\ast} + \Delta t_j^z,
\end{equation}
or 
\begin{equation} \label{eq:B_iz2}
t_i^{z^\ast} + \Delta t_i^z \leq t_j^{z^\ast}.
\end{equation}

Let CAV $k\in\mathcal{A}_i^l,~l=l_i^{f^\ast}$, be the vehicle immediately ahead of CAV $i\in\mathcal{N}(t)$ at lane $l_i^{f^\ast}\in\mathcal{L}$. In order for CAV $i$ to avoid the rear-end collision at the merging zone $z\in\mathcal{Z}_i$,
\begin{equation}\label{eq:rearEndatmz}
    {t_{i}^{z^\ast}}\geq{t_{k}^{z^\ast}} + \rho_k^z, 
\end{equation}
where $\rho_k^z\in\mathbb{R}^+$ is the time that it takes for CAV $k$ to travel a safe distance $\delta$ inside the merging zone $z$.

As we mentioned earlier, in the upper-level planning, we relax the FIFO queuing policy to improve the traffic throughput in multiple intersections. Upon arrival at the control zone, each CAV $i\in\mathcal{N}(t)$ recursively computes the energy-optimal arrival time for each merging zone along its path ensuring the lateral safety in conjunction with the lane that it should occupy using the Algorithms \ref{Alg:Scheduling} and \ref{Alg:lane-changing}. Given $l_i^f$, the lane that CAV $i$ needs to follow after lane changing zone, CAV $i$ employs Algorithm \ref{Alg:Scheduling} to find the energy-optimal arrival time for all merging zones $z\in\mathcal{Z}_i$.

\begin{algorithm}
 \caption{Arrival time at merging zones}
\hspace*{\algorithmicindent} \textbf{Input:} $l_i^f$ \\
\hspace*{\algorithmicindent} \textbf{Output:} Arrival time at merging zones $\mathcal{Z}_i = \{z_1,\dots,z_n\}$ 
%\renewcommand{\algorithmicrequire}{\textbf{Input:}}
 %\renewcommand{\algorithmicensure}{\textbf{Output:}}
 %\REQUIRE $l_i^f$
 %\ENSURE  arrival time at merging zones $\mathcal{Z}_i = \{z_1,\dots,z_n\}$
 \begin{algorithmic}[1]

 \For{ $z\in\mathcal{Z}_i$}
 \State $k=\max\{j~|~j\in\mathcal{A}_i^{l}~,~l = l_i^f\}$ \label{alg:check k}
 \State $t_i^{z^\ast}\gets \max\left\{\bar{t}_i^{z}, {t_{k}^{z^\ast}} + \rho_k^z \right\}$\label{alg:tbar,k}
 \State Sort $\mathcal{T}_i^z$ increasingly
  \For {$t_j^{z^\ast}\in \mathcal{T}_i^z$} \label{alg:forMZStart}
    \If {$t_j^{z^\ast}+\Delta t _j^z  \leq t_i^{z^\ast}$}\label{alg:tbarcont1}
  \State $continue$ \Comment{CAV $j$ exits before CAV $i$ enters}
  \EndIf\label{alg:tbarcont1End}
  \If {$t_i^{z^\ast}+\Delta t _i^z  \leq t_j^{z^\ast}$}
  \State $break$ \Comment{CAV $i$ exits before CAV $j$ enters}
  \Else\Comment{CAV $i$ has conflict with CAV $j$}
  \State 	$t_i^{z^\ast}\gets t_j^{z^\ast}+\Delta t _j^z $
  \EndIf
  \EndFor \label{alg:forMZEND}
  \EndFor
 \State \Return $\{t_i^{z^\ast}~|~z\in\mathcal{Z}_i\}$
 \end{algorithmic} \label{Alg:Scheduling}
 \end{algorithm}

\begin{theorem}
For a given lane $l_i^f\in\mathcal{L}$, CAV $i\in\mathcal{N}(t)$ recursively computes the energy-optimal arrival time at merging zone $z \in\mathcal{Z}_i$ subject to the constraints \eqref{eq:B_iz}-\eqref{eq:rearEndatmz} using Algorithm \ref{Alg:Scheduling}. 
\end{theorem}
\begin{proof}
For a given lane $l_i^f$ after lane-changing zone, there are four cases to consider.

Case 1: If $\mathcal{N}(t_i^0)$ = $\mathcal{C}_i$, then we have: $\mathcal{B}_i^z=\emptyset$ for all $z\in\mathcal{Z}_i$ and $\mathcal{A}_i^l=\emptyset$ for all $l\in\mathcal{L}$. Thus, there is no leading vehicle (in line \ref{alg:check k}, $k=\emptyset$), and $ \mathcal{T}^z_i=\emptyset$. Thus, we have $t_i^{z^\ast}=\bar{t}_i^{z}$ for all $z\in\mathcal{Z}_i$, and from Lemma \ref{lem:tbar}, $\bar{t}_i^{z}$ is the energy-optimal arrival time.\\
Case 2: If $\mathcal{N}(t_i^0)=\mathcal{A}_i^l,~l=l_i^f\in\mathcal{L}$, there exists a CAV $k\in\mathcal{A}_i^{l},~l=l_i^f$, which is immediately ahead of CAV $i$ at lane $l_i^f$ (in line \ref{alg:check k}, $k\neq\emptyset$). Similarly, we have: $\mathcal{B}_i^z=\emptyset$ for all $z\in\mathcal{Z}_i$ and $\mathcal{C}=\emptyset$. To ensure rear-end safety, we have: 
\begin{equation}
    \max\left\{\bar{t}_i^{z}, {t_{k}^{z^\ast}} + \rho_k^z\right\}\leq t_i^{z^\ast},~\forall z\in\mathcal{Z}_i.
\end{equation}
Selecting the lower bound in the above equation, CAV $i$ computes the energy-optimal arrival time at the merging zone $z\in\mathcal{Z}_i$ satisfying the rear-end safety constraint \eqref{eq:rearEndatmz} (line \ref{alg:tbar,k}).
\\
Case 3: If $\mathcal{N}(t_i^0)=\bigcup\limits_{z\in\mathcal{Z}_i}\mathcal{B}_i^z$, then we have: $\mathcal{C}=\emptyset$ and $\mathcal{A}_i^l=\emptyset$ for all $l\in\mathcal{L}$.  Since for CAV $i,$ $\bar{t}_i^{z}$ is the energy-optimal arrival time at merging zone $z\in\mathcal{Z}_i$ without considering safety (Lemma \ref{lem:tbar}), we use $\bar{t}_i^{z}$ as a lower bound (line \ref{alg:tbarcont1}-\ref{alg:tbarcont1End}). CAV $i$ determines $t_i^{z^\ast}\in[\bar{t}_i^z,\infty]$ to be the smallest time satisfying lateral safety constraints in \eqref{eq:B_iz} or \eqref{eq:B_iz2} (line \ref{alg:forMZStart} -\ref{alg:forMZEND}).
\\
Case 4: Similarly, for other cases that $\mathcal{A}_i^l\neq \emptyset,~l\in\mathcal{L}$, $\mathcal{B}_i^z\neq \emptyset~z\in\mathcal{Z}_i$, and $\mathcal{C}_i\neq \emptyset$, the energy-optimal arrival time at the merging zone $z\in\mathcal{Z}_i$, $t_i^{z^\ast}\in[\bar{t}_i^z,\infty]$ is the smallest time satisfying the lateral safety constraints in \eqref{eq:B_iz} or \eqref{eq:B_iz2} (line \ref{alg:forMZStart} -\ref{alg:forMZEND}), and rear-end safety constraint \eqref{eq:rearEndatmz} (line \ref{alg:tbar,k}).
\end{proof} 

\begin{algorithm}
 \caption{Lane-changing decision}
 \begin{algorithmic}[1]
 \State occupancy $\gets \bigcup \limits_{j\in\mathcal{A}_i^l,l\in\mathcal{L}}[t_i^0, \infty) \cap\Gamma_j$\label{al:occu}
 \State $l_i^{f^\ast}\gets l_i^0$
 \State $t_i^{z_n^{\ast}}\gets$ Find arrival time at $z_n$, given $l_i^{f^\ast}$
 \If{occupancy $=\emptyset$}
 \For{$l\in\mathcal{L}\setminus {l_i^0}$ }
 \State $l_i^f\gets l$
 \State $t_i^{z_n}\gets$ Find arrival time at $z_n$, given $l_i^f$
  \If{$t_i^{z_n}<t_i^{z_n^{\ast}}$}
  \State $l_i^{f^\ast}\gets l$
  \State $t_i^{z_n^{\ast}}\gets t_i^{z_n}$
  \EndIf
  \EndFor
  \EndIf
 \State return $l_i^{f^\ast}$
 \end{algorithmic}\label{Alg:lane-changing}
 \end{algorithm}

Using  Algorithm \ref{Alg:lane-changing}, CAV $i\in\mathcal{N}(t)$ investigates the feasibility of the lane-changing maneuver (line \ref{al:occu}). If such a maneuver is not feasible, CAV $i$ exits from the lane-changing zone with the same lane that it entered the control zone, $l_i^{f^\ast}=l_i^0\in\mathcal{L}$, and no lane-change maneuver is performed. Otherwise, CAV $i$ selects the optimal lane $l_i^{f^\ast}\in\mathcal{L}$ yielding the minimum travel time, i.e., the arrival time at the last merging zone.

\section{Low-level Planning}\label{section4}
\begin{figure}[ht]
\centering
\includegraphics[width=0.8\linewidth]{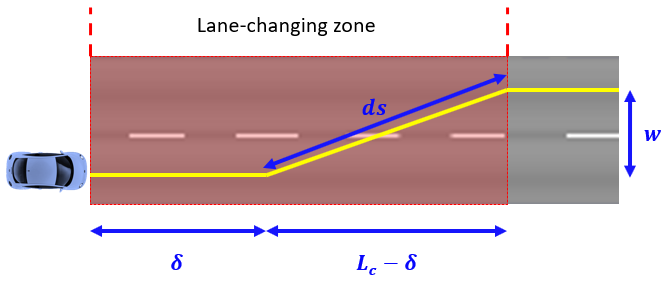}
\caption{Lane-changing maneuver. }
\label{fig:lanechangin}%
\end{figure}

In our decentralized framework, the outputs of the upper-level planning for CAV $i\in\mathcal{N}(t)$, which are the optimal arrival time, $t_i^{z^\ast},$ at merging zone $z\in\mathcal{Z}_i$ along with the optimal lane to occupy, $l_i^{f^\ast}$, become inputs for the low-level planning. In particular, in the low-level planning, CAV $i\in\mathcal{N}(t)$ formulates the optimal control problem with interior-point constraints defined at the entry and exit of each merging zone, the solution of which minimizes the engine control effort, and energy consumption correspondingly. 
In addition, if the optimal lane after the lane-changing zone, $l_i^{f^\ast}\neq l_i^0$, CAV $i$ needs to perform the lane-changing maneuver. If so, CAV $i\in\mathcal{N}(t)$ first travels a distance $\delta\in\mathbb{R}^+$ on the lane that it entered, $l_i^0\in\mathcal{L}$, and then travels on a triangle with the hypotenuse $ds\in\mathbb{R}^+$ and sides $Lc - \delta$ and $w\in\mathbb{R}^+$, respectively, to reach lane $l_i^{f^\ast}$ (see Fig. \ref{fig:lanechangin}). Since $w\ll L_c-\delta$, we approximate $ds\approx L_c-\delta$.

For CAV $i\in\mathcal{N}(t)$ the control-effort minimization with interior-point constraints at the boundary of each merging zone $z\in\mathcal{Z}_i$ is formulated as follows:
\begin{problem}\label{problem1} Control-effort minimization
\begin{equation}\label{EnergyOptimalProblem}
 \begin{array}{ll}
 \min\limits_{\textit{u}_i\in\mathcal{U}_i} \quad &J_i(u_i(t))= {\dfrac{1}{2}} \bigintss_{t_i^{0}}^{t_i^f} u_{i}(t)^2dt, \\
\emph{subject to:}&\quad(\ref{27a}), (\ref{uconstraint}), (\ref{vconstraint}),(\ref{RearEndCons}),\\
\emph{given }&\quad p_i(t_i^{0}), v_i(t_i^{0}),\\
&\quad p_i(t_i^{z^\ast}),p_i(t_i^{z^\ast}+\Delta t_i^{z}), ~~\forall z\in\mathcal{Z}_i.\\
           \end{array} 
\end{equation}\end{problem}
Recall that $t_i^f$ is the time that CAV $i$ exits the control zone, i.e., the merging zone of the last intersection along its path, $t_i^f = t_i^{z_n^\ast} + \Delta t_i^{z_n}$.

\subsection{Solution of the Control-Effort Minimization}
To derive the solution of the control-effort minimization (Problem \ref{problem1}), we apply Hamiltonian analysis. After solving the upper-level problem for CAV $i\in\mathcal{N}(t)$, the entry and exit time of merging zones $z\in\mathcal{Z}_i$ are the interior-point constraints for the low-level problem. First, we adjoin the control inequality constraints \eqref{uconstraint} along with the $q$th-order state variable inequality constraints to the Hamiltonian function. 
The $q$th-order state variable inequality constraint can be found by taking the successive total time derivative of constraint and substitute \eqref{27a}, until we obtain an expression that is explicitly dependent on the control variable \cite{bryson1975applied}. 
For each CAV $i\in\mathcal{N}(t)$, with CAV $k\in\mathcal{A}_i^l,~l=l_i^{f^\ast}$, physically located ahead of it, the Hamiltonian given by
\begin{equation}\label{c1}
\begin{aligned}
H_i&(t,p_i(t),v_i(t),u_i(t))=\frac{1}{2}u_i(t)^2+\lambda_i^p v_i (t)+\lambda_i^v u_i(t)\\
+&\mu_i^a(u_i (t)-u_{{i,\max}})+\mu_i^b(u_{i,\min}-u_i(t))\\
+&\mu_i^c(u_i(t))+\mu_i^d(-u_i(t))\\
+&\mu_i^s(u_i(t)-u_k^\ast(t)),
\end{aligned}
\end{equation}
where \(\lambda_i^p\) and \(\lambda_i^v\) are costates, and $\mathbf{\mu}_i^\top = [\mu_i^a,\mu_i^b,\mu_i^c,\mu_i^d,\mu_i^s]$ is a vector of Lagrange multipliers. It should be noted that $u_k^\ast(t)$ is the optimal control input for CAV $k\in\mathcal{A}_i^l,~l=l_i^{f^\ast}$, which is available information to CAV $i$ through the coordinator.
\\
The Euler-Lagrange equations become: 
\begin{align}
\dot{\lambda}_i^p&=-\frac{\partial H_i}{\partial p_i}=0,\label{euler:optimallambdap}\\
\dot{\lambda}_i^v&=-\frac{\partial H_i}{\partial v_i}= -\lambda_i^p,\label{euler:optimallambdav}\\
\frac{\partial H_i}{\partial u_i}&=u_i+\lambda_i^v+\mu_i^a-\mu_i^b+\mu_i^c-\mu_i^d+\mu_i^s=0\label{euler:optimalU}.
\end{align}
Since the speed of CAV $i$ is not specified at the fixed terminal time $t_i^f$, we have \cite{bryson1975applied}: 
\begin{equation}
    \lambda_i^v(t_i^f) = 0.
\end{equation}

Since the arrival time of CAV $i\in\mathcal{N}(t)$ at the entry and exit of each merging zone $z\in\mathcal{Z}_i$, is specified, for each interior-point constraint at specified time $t_1$, we have the following condition
\begin{equation}
    \mathbf{N}(\mathbf{x}_i(t),t) = \left[ \begin{array}{c}
        p_i(t) - C\\
        t - t_1
    \end{array}
\right]=0
    ,\label{eq:Interior}
\end{equation}
where $C$ is the position at the specified time $t_1$ (i.e., entry/exit of merging zone). In addition, the costates and the Hamiltonian should satisfy the following jump conditions at $t_1^-$ and $t_1^+$, \cite{bryson1975applied}

\begin{align}
\mathbf{\lambda}_i^\top(t_1^-) &= \mathbf{\lambda}_i^\top(t_1^+) + \mathbf{\pi}^\top~ \frac{\partial \mathbf{N}}{\partial \mathbf{x}_i}|_{t=t_1},\\
    H_i(t_1^-)&=H_i(t_1^+)-\pi^\top~\frac{\partial \mathbf{N}}{\partial t}|_{t=t_1}. 
\end{align}
where $\lambda_i^\top = [\lambda_i^p, \lambda_i^v]$, $\pi^\top = [\pi_1, \pi_2]$, $\dfrac{\partial \mathbf{N}}{\partial x_i} = \left[\begin{array}{cc}
     1&0  \\
     0&0 
\end{array}\right]$ and $\dfrac{\partial \mathbf{N}}{\partial t} = [0,1]^\top$. Hence, 
\begin{align}
    \lambda_i^p(t_1^-) &=\lambda_i^p(t_1^+)+ \pi_1\label{lambdap=},\\
    \lambda_i^v(t_1^-)&=\lambda_i^v(t_1^+)\label{lambdav=},\\
    H_i(t_1^-)&=H_i(t_1^+)-\pi_2\label{Hamil}. 
\end{align}
Note that $\pi^\top$ is a 2-component vector of constant Lagrange multipliers, determined so that the interior-point constraint \eqref{eq:Interior} is satisfied.
\subsubsection{Unconstrained Solution Without Interior-point Constraints}
If the state and control constraints never become active, \(\mu_i^a=\mu_i^b=\mu_i^c=\mu_i^d=\mu_i^s=0\) the solution, see \cite{Malikopoulos2017}, is
\begin{equation}\label{27}
 u_i^{\ast}(t)=a_it+b_i,
\end{equation}
by substituting \eqref{27} in \eqref{27a}, we have
\begin{align}\label{27b}
&v_i^{\ast}(t)=\frac{1}{2}a_it^2+b_it+c_i,\\
&p_i^{\ast}(t)=\frac{1}{6}a_it^3+\frac{1}{2}b_it^2+c_it+d_i.
\end{align}
In the above equations \(a_i,b_i,c_i,d_i\) are constants of integration, which are found by substituting the initial and final conditions $p^\ast_i(t_i^0), v^\ast_i(t_i^0)$, $p^\ast_i(t_i^f)$ and $u^\ast_i(t_i^f) = 0$.

\subsubsection{Unconstrained Solution With Interior-point Constraints}
To find the analytical solution for CAV $i\in\mathcal{N}(t)$ including the interior-point constraints at the entry and exit of merging zone $z\in\mathcal{Z}_i$ (recall that $z_1$ and $z_n$ are the first and last merging zones that CAV $i$ crosses, respectively, Definition \ref{Defn:Z_i}), we need to satisfy $2n-1$ interior-point constraints ($t_i^f$ is excluded, since it is a boundary condition). 
\begin{lemma}\label{lem: continuity}
The optimal control input $u_i^\ast$(t) when none of the constraints is active at the interior-point constraint $N_j,~j\in\{1,\dots,2n-1\}$, where $n$ is the total number of merging zones in CAV $i$'s path, is continuous.
\end{lemma}
\begin{proof}
Let $t_j$ be the time that we have an interior-point constraint $N_j,~j\in\{1,\dots,2n-1\}$. 
From \eqref{lambdav=}, we know $\lambda^v_i$ is continuous $t_j$, i.e.,
\begin{equation}
    \lambda_i^v(t_j^-) =  \lambda_i^v(t_j^+),\\
\end{equation}
 Since none of the state or control inequality constraints is active, we have \(\mu_i^a=\mu_i^b=\mu_i^c=\mu_i^d=\mu_i^s=0\) and from \eqref{euler:optimalU}, we have $\lambda_i^v(t) = - u_i(t)$ which gives 
\begin{equation}
    u_i(t_j^-) = u_i(t_j^+).
\end{equation}
\end{proof}
\begin{theorem}\label{theorem2}
The unconstrained solution of Problem \ref{problem1} for CAV $i\in\mathcal{N}(t)$ with $n$ merging zones, is a continuous piecewise linear function. 
\end{theorem}
\begin{proof}
For CAV $i\in\mathcal{N}(t)$, we first divide $[t_i^0,t_i^f]$ (recall that $t_i^f=t_i^{z_n^\ast} + \Delta t_i^{z_n}$) into $2n$ sub-intervals with corresponding optimal control input as follows: 
\begin{equation}\label{eq:Utotal0}
 u^*_i(t)=\left\{ \begin{array}{ll}
u^{(1)}_i(t),& \mbox{if}\quad t_i^0\leq t<t_i^{z_1^\ast},\\
u^{(2)}_i(t), & \mbox{if}\quad t_i^{z_1^\ast}\leq t\leq t_i^{z_1^\ast}+\Delta t_i^{z_1},\\
   &\vdots\\
u^{(2n-1)}_i(t), & \mbox{if}\quad t_i^{z_{n-1}^\ast}+ \Delta t_i^{z_{n-1}}\leq t\leq t_i^{z_{n}^\ast},\\
u^{(2n)}_i(t), & \mbox{if}\quad t_i^{z_{n}^\ast}\leq t\leq t_i^f.\\
           \end{array} \right. 
\end{equation}
Integrating \eqref{euler:optimallambdap} and \eqref{euler:optimallambdav} at each time-interval $j\in\{1,\dots,2n\}$ and using \eqref{euler:optimalU}, we get a linear form for $u^{(j)}_i(t)$. From Lemma \ref{lem: continuity}, we have continuity of control input at each interior point, thus the control input is a continuous piecewise linear function.
\end{proof} 

\begin{corollary}\label{cor= pi2}
For CAV $i\in\mathcal{N}(t)$, let ${\pi^{(j)}} = [\pi^{(j)}_1, \pi^{(j)}_2]^\top$ be the constant Lagrange multipliers for the interior-point constraint $N_j,~j\in\{1,\dots,2n-1\}$ at time $t_j$, where $n$ is the total number of merging zones in CAV $i$'s path. Then, we have
\begin{equation}
    \pi^{(j)}_2= -\pi^{(j)}_1 v_i(t_j) .
\end{equation}
\end{corollary}
\begin{proof}
From \eqref{Hamil} and \eqref{c1}, we have 
\begin{align}\label{NJump}
    \frac{1}{2} u_i(t_j^-)^2 + \lambda_i^p(t_j^-) v_i(t_j^-) + \lambda_i^v(t_j^-) u_i(t_j^-) =\nonumber\\ \frac{1}{2} u_i(t_j^+)^2 + \lambda_i^p(t_j^+) v_i(t_j^+) + \lambda_i^v(t_j^+) u_i(t_j^+)-\pi_2,
\end{align}
and by substituting \eqref{lambdap=} and $\lambda_i^v(t) = - u_i(t)$ into \eqref{NJump}, we get
\begin{align}\label{Njump2}
    -\frac{1}{2} u_i(t_j^-)^2 + (\lambda_i^p(t_j^+)+\pi^{(j)}_1)~v_i(t_j^-)  =\nonumber\\ -\frac{1}{2} u_i(t_j^+)^2 + \lambda_i^p(t_j^+) v_i(t_j^+) -\pi^{(j)}_2.
\end{align}
Using continuity of speed at the interior-point, i.e., $v_i(t_j^-)=v_i(t_j^+) = v_i(t_j)$ and rearranging \eqref{Njump2}, we get
\begin{equation}
    -\frac{1}{2} u_i(t_j^-)^2 + \pi^{(j)}_1 v_i(t_j) = -\frac{1}{2} u_i(t_j^+)^2  -\pi^{(j)}_2.
\end{equation}
From Lemma \eqref{lem: continuity}, we have $u_i(t_j^-) = u_i(t_j^+)$. Therefore $\pi^{(j)}_2= -\pi^{(j)}_1 v_i(t_j)$, and the proof is complete.
\end{proof}

The unconstrained solution with interior-point constraints for CAV $i\in\mathcal{N}(t)$ with $n$ merging zones, consists of $2n$ unconstrained arcs as follows: 
\begin{equation}\label{eq:Utotal}\small{
 u^*_i(t)=\left\{ \begin{array}{ll}

a_i^{(1)}t+b_i^{(1)},& \mbox{if}\quad t_i^0\leq t<t_i^{z_1^\ast},\\
a_i^{(2)}t+b_i^{(2)},& \mbox{if}\quad t_i^{z_1^\ast}\leq t\leq t_i^{z_1^\ast}+\Delta t_i^{z_1},\\
   &\vdots\\
a_i^{(2n-1)}t+b_i^{(2n-1)}, & \mbox{if}\quad t_i^{z_{n-1}^\ast}+ \Delta t_i^{z_{n-1}}\leq t\leq t_i^{z_{n}^\ast},\\
a_i^{(2n)}t+b_i^{(2n)}, & \mbox{if}\quad t_i^{z_{n}^\ast}\leq t\leq t_i^f,\\
           \end{array} \right. }
\end{equation}
where $a_i^{(j)}$ and $b_i^{(j)}$ are unknown parameters for the unconstrained arc $j\in\{1,\dots,2n\}$. Substituting \eqref{eq:Utotal} in \eqref{27a}, and integrating, we get two more unknowns per arc that are constants of integration. Therefore, for $n$ merging zones we have $4(2n)=8n$ unknowns that need to be computed along with $2n-1$ constant Lagrange multipliers ($\pi^{(j)}_1, j\in\{1,\dots,2n-1\}$) resulting in total $10n-1$ unknowns. 
Initial conditions $p_i(t_i^0)$ and $v_i(t_i^0)$, final conditions $p_i(t_i^f)$ and $u_i(t_i^f)$, (4 equations), continuity of state and control at interior-point constraints ($3\cdot(2n-1)$ equations), position at interior-point constraints ($2n-1$ equations), and jump conditions on $\lambda_i^p$ at interior-point constraints ($2n-1$ equations) result in $10n-1$ equations, which form a system of linear equations that yields the optimal trajectory.
\begin{remark}
After finding the optimal trajectory, along with $\pi^{(j)}_1, j\in\{1,\dots,2n-1\}$, one can use Corollary \ref{cor= pi2} to find $\pi^{(j)}_2$ which satisfies \eqref{Hamil}.
\end{remark}

\begin{theorem}
For the cases that none of the state/control inequality constraints becomes active, the control effort minimization problem with interior-point constraints always has a unique solution. 
\end{theorem}
\begin{proof}
The analytical solution for the unconstrained case with interior-point constraints is found by solving a system of linear equations in the classical form $AX=B$, where $A\in\mathbb{R}^{(10n-1)\times (10n-1)}$ is the coefficients matrix, $X\in\mathbb{R}^{(10n-1)\times 1}$ is the vector of $10n-1$ unknowns and $b\in\mathbb{R}^{(10n-1)\times 1}$ is the constant vector. Since there are $10n-1$ linearly independent equations, which forms row vectors in $AX=B$, we have $\text{rank}(A) = \text{rank}(A|B) = 10n-1$, and the proof is complete.
\end{proof}
\subsubsection{Constrained solution}
Using \eqref{eq:Utotal}, we first start with the unconstrained solution of Problem \ref{problem1}. If the solution violates any of the speed \eqref{vconstraint} or control \eqref{uconstraint} constraints, then the unconstrained arc is pieced together with the arc corresponding to the violated constraint at unknown time $\tau_1$, and we re-solve the problem with the two arcs pieced together. The two arcs yield a set of algebraic equations which are solved simultaneously using the boundary and interior conditions at $\tau_1$. If the resulting solution violates another constraint, then the last two arcs are pieced together with the arc corresponding to the new violated constraint, and we re-solve the problem with the three arcs pieced together at unknown times $\tau_1$ and $\tau_2$. The three arcs will yield a new set of algebraic equations that need to be solved simultaneously using the boundary and interior conditions at $\tau_1$ and $\tau_2$. The process is repeated until the solution does not violate any other constraints, \cite{Malikopoulos2017,chalaki2020TITS}. In the following section, we show the analysis for the case where the rear-end safety constraint becomes active.
\subsubsection{Rear-end safety constraint becomes active}
Suppose for CAV $i\in\mathcal{N}(t)$, at some time $t=\tau_1$, the rear-end safety constraint with the vehicle $k$ becomes active until $t =\tau_2$, $p_k(t)-p_i(t) = \delta$ for all $t\in[\tau_1,\tau_2]$, in this case $\mu_i^s\neq0$.
At the entry of the constrained arc, we have the following tangency conditions

\begin{equation}
    \mathbf{N}(\mathbf{x}_i(t),t) = \left[ \begin{array}{c}
        p_i(t)-p_k^\ast(t)+\delta\\
        v_i(t)-v_k^\ast(t)
    \end{array}
\right]=0
    .\label{eq:tangency}
\end{equation}
 Since $\mathbf{N}(t,\mathbf{x}_i(t)) = 0$ for $t\in[\tau_1,\tau_2]$, its first derivative, which is dependent on the optimal control input, should vanish in $t\in[\tau_1,\tau_2]$, i.e.,
\begin{equation}\label{1-storderSafety}
    \mathbf{N}^{(1)}(t,\mathbf{x}_i(t))=u^{\ast}_i(t)-u^{\ast}_k(t)= 0.
\end{equation}
From \eqref{1-storderSafety}, the optimal control input of CAV $i\in\mathcal{N}(t)$, when rear-end safety constraint is active, can be found $u^{\ast}_i(t)=u^{\ast}_k(t)$
The optimal solution needs to satisfy the following jump conditions on costates upon entry to the constrained arc at $t=\tau_1$,
\begin{align}
\lambda_i^p(\tau_1^-) &= \lambda_i^p(\tau_1^+)+[\pi_1,\pi_2]\frac{\partial \mathbf{N}}{\partial p_i} = \lambda_i^p(\tau_1^+)+\pi_1,\\
\lambda_i^v(\tau_1^-) &= \lambda_i^v(\tau_1^+)+[\pi_1,\pi_2]\frac{\partial \mathbf{N}}{\partial v_i} = \lambda_i^v(\tau_1^+)+\pi_2,\\
H_i(\tau_1^-) &= H_i(\tau_1^+)-[\pi_1,\pi_2]\frac{\partial \mathbf{N}}{\partial t} \nonumber\\&= H_i(\tau_1^+)+\pi_1v^{\ast}_k(t)+\pi_2u^{\ast}_k(t),\label{37safety}
\end{align}
where $\pi_1$ and $\pi_2$ are constant Lagrange multipliers, determined so that \eqref{eq:tangency} is satisfied.
At the exit point of the constrained arc, we have
\begin{align}
\lambda_i^p(\tau_2^-) &= \lambda_i^p(\tau_2^+),\label{39safety}\\
\lambda_i^v(\tau_2^-) &= \lambda_i^v(\tau_2^+),\label{40safety}\\
H_i(\tau_2^-) &= H_i(\tau_2^+).\label{41safety}
\end{align}
As described earlier, the three arcs need to be solved simultaneously using initial and final conditions (speed and position), interior-point constraints \eqref{eq:Interior}, and interior conditions at unknown times $\tau_1$ and $\tau_2$ (continuity of speed and position, jump conditions \eqref{eq:tangency}-\eqref{41safety}).

\section{Deviation From Nominal Planned Position}\label{section5}
In this section, we extend the previous results to include a bounded steady-state error in CAV's position, which can be originated from the vehicle-level controller tracking the optimal trajectory. Namely, suppose that CAV $i$'s actual position deviates from the nominal $p^\ast_i(t)$, which is the optimal solution of the Problem \ref{problem1}, and it takes values in $[p^\ast_i(t)-\epsilon~,p^\ast_i(t)+\epsilon]$, where $\epsilon\in\mathbb{R}^+$ is the maximum deviation from the nominal path. To guarantee longitudinal and lateral safety, we consider the worst-case scenario in our upper-level and low-level planning analysis.
\subsection{Low-level Safety}

To guarantee rear-end safety between CAV $i\in\mathcal{N}(t)$ and CAV $k\in\mathcal{N}(t)$, where CAV $k\in\mathcal{A}_i^l,~l=l_i^{f^\ast}$ is the vehicle immediately ahead of CAV $i$, we modify the rear-end safety constraint \eqref{RearEndCons} as follows: 

\begin{equation}
    (p_k(t)-\epsilon) - (p_i(t)+\epsilon) \geq \delta,
\end{equation}
where it simplifies to 
\begin{equation}\label{newrear}
    p_k(t) - p_i(t) \geq \delta + 2\epsilon.
\end{equation}
Thus, considering the worst-case scenario in the low-level planning results in increasing the rear-end safety distance $\delta_{\text{new}}=\delta + 2\epsilon$. However, this might be a conservative causing the rear-end safety constraint becomes active without being necessary and thus results in higher fuel consumption. 
\subsection{Upper-level Safety} 
In the low-level problem, we modified the rear-end safety constraint to guarantee safety in the worst-case scenario in the presence of a bounded steady-state error in the position. In addition, to ensure safety in the upper-level problem, we need to consider the error's effects in order to avoid lateral collision. By introducing the idle-time $t_{\text{idle}}$, acting as a safety buffer in the presence of a bounded steady-state error in position, we modify the lateral safety constraints with idle-time as follows   
\begin{equation} \label{eq:B_izIDLE1}
t_i^{z^\ast} \geq t_j^{z^\ast} + \Delta t_j^z + t_\text{idle},
\end{equation}
or 
\begin{equation} \label{eq:B_izIDLE2}
t_i^{z^\ast} + \Delta t_i^z + t_\text{idle} \leq t_j^{z^\ast},
\end{equation}
where in \eqref{eq:B_izIDLE1}, merging zone should be idle after CAV $j$'s planned exit time and in \eqref{eq:B_izIDLE2} merging zone should be idle before CAV  $j$'s planned arrival time. 
To consider rear-end safety at entry of each merging zone, \eqref{eq:rearEndatmz} is being adjusted as follows  
\begin{equation}\label{eq:rearEndatmzIDLE}
    {t_{i}^{z^\ast}}\geq{t_{k}^{z^\ast}} + \rho_k^z + t_{\text{idle}}.
\end{equation}

The worst-case scenario for computing $t_\text{idle}$ can be computed as follows 
\begin{equation}\label{t_idle}
    t_\text{idle} = \frac{2\epsilon}{v_\text{min}}.
\end{equation}
It should be noted that $t_{\text{idle}}$ computed from \eqref{t_idle} is very conservative, since it assumes that both vehicles cross the merging zone with minimum speed. 
Algorithm \ref{Alg:idleTIME} can be employed to consider the deviation from the nominal planned position. 

\begin{algorithm}
 \caption{Arrival time at merging zones with idle-time}
\hspace*{\algorithmicindent} \textbf{Input:} $l_i^f$\\ 
\hspace*{\algorithmicindent} \textbf{Output:} Arrival time at merging zones $\mathcal{Z}_i = \{z_1,\dots,z_n\}$ 
 \begin{algorithmic}[1]
 \For{ $z\in\mathcal{Z}_i$ }
 \State $k = \max\{j~|~j\in\mathcal{A}_i^{l}~,~l = l_i^f\}$
 \State $t_i^{z^\ast}\gets \max\left\{\bar{t}_i^{z}, {t_{k}^{z^\ast}} + \rho_k^z+t_{\text{idle}} \right\}$
  \State Sort $\mathcal{T}_i^z$ increasingly
  \For {$t_j^{z^\ast}\in \mathcal{T}_i^z$} 
  \If {$t_j^{z^\ast}+\Delta t _j^z + t_{\text{idle}} \leq t_i^{z^\ast}$}
  \State $continue$
  \EndIf
  \If {$t_i^{z^\ast}+\Delta t _i^z + t_{\text{idle}} \leq t_j^{z^\ast}$}
  \State $break$
  \Else
  \State 	$t_i^{z^\ast}\gets t_j^{z^\ast}+\Delta t _j^z+t_{\text{idle}} $;
  \EndIf
  \EndFor
  \EndFor
 \State return $\{t_i^{z^\ast}~|~z\in\mathcal{Z}_i\}$
 \end{algorithmic} \label{Alg:idleTIME}
 \end{algorithm}
 
\section{Simulation Example}\label{section6}
To evaluate the effectiveness of the proposed framework in reducing fuel consumption and improving traffic throughput, we investigate the coordination of CAVs at three adjacent intersections in two scenarios under different traffic volumes, and then compare the results with a baseline scenario consisting of two-phase traffic signals.
We construct the baseline scenario with two-phase fixed-time traffic signals in PTV-VISSIM \cite{PTV}, which is a commercial microscopic multi-modal traffic flow simulation software, by considering all vehicles as human-driven and without V2V communication. To emulate human-driven vehicles' driving behavior, we use a built-in car-following model (Wiedemann \cite{Wiedemann1974}) in PTV-VISSIM with default parameters. In the optimal-scenario, we use a dynamic-link library in PTV-VISSIM to simulate our framework. Videos from our simulation analysis can be found at the supplemental site, \url{https://sites.google.com/view/ud-ids-lab/OCMI}.

\begin{figure}[ht]
\centering
\includegraphics[width=0.95\linewidth]{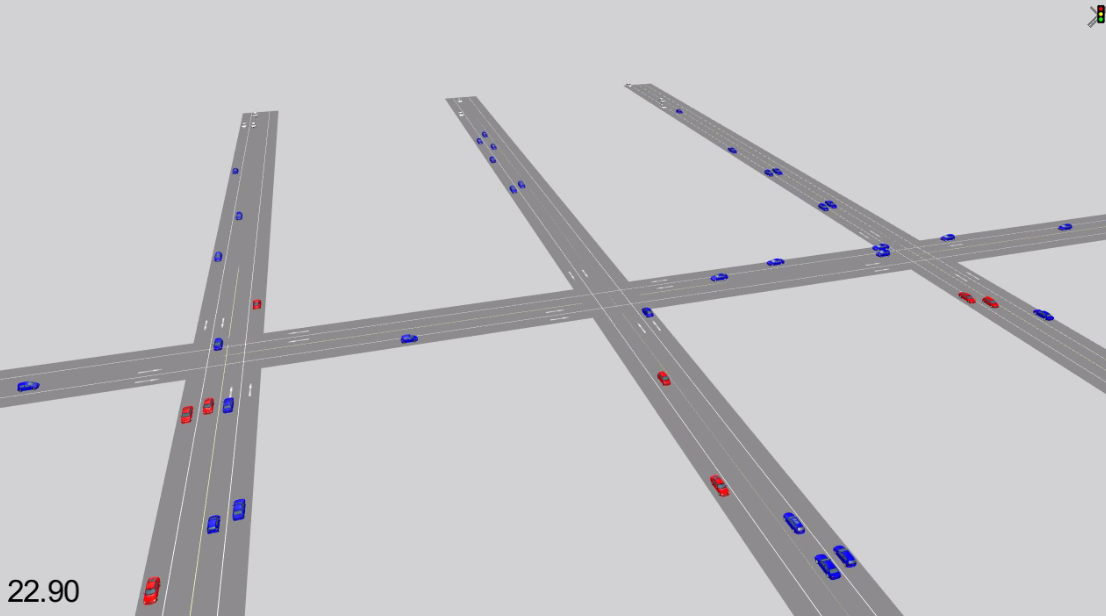}
\caption{Snap-shot of the three multi-lane adjacent intersections in the first scenario.}
\label{fig:Simulation}%
\end{figure}

For the first scenario, we construct three symmetric adjacent intersections. We consider the length of each road connecting to the intersections to be $L=150$ m, the length of the merging zones to be $S=15$ m, and the distance between each intersection to be $D=75$ m (see Fig. \ref{fig:Simulation}). The CAVs enter the control zone with initial speed uniformly distributed between $11$ m/s to $13$ m/s from each entry with equal traffic volumes. 

\begin{figure}[ht]
\centering
\includegraphics[width=0.95\linewidth]{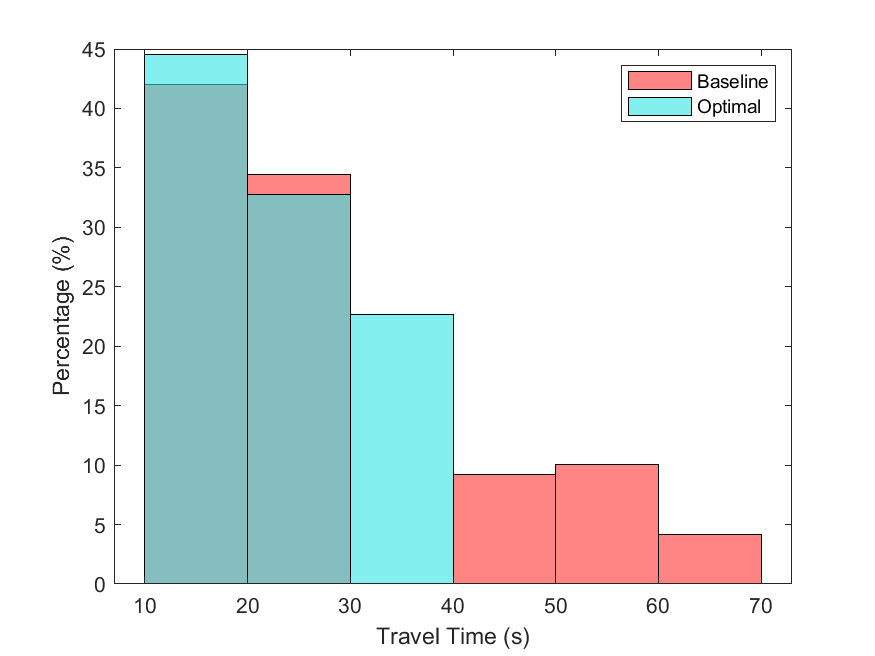}
\caption{A relative frequency histogram for travel time of each vehicle for the first scenario in the baseline and optimal cases under traffic volume $\num{1400}$ veh/h. }
\label{fig:histTravelTime}%
\end{figure}

Table \ref{tbl:average travel time} shows the average travel time of all CAVs inside the control zone for the baseline and optimal cases, respectively, at different traffic volumes ranging from $600$ veh/h to $\num{1400}$ veh/h per lane for each entry. For each traffic volume, we performed five simulations with different random seeds and averaged the results. Within our proposed framework, average travel time has been decreased by $11\% -24\%$ compared to the baseline scenario with two-phase traffic signals. Relative frequency histogram of travel time of each CAV for traffic volume $\num{1400}$ veh/h for one of the selected seed, which includes $119$ vehicles, for the baseline and optimal scenarios are shown in Fig. \ref{fig:histTravelTime}. 
As it can be seen in Fig. \ref{fig:histTravelTime}, travel times for all CAVs are less than $40$ s for the optimal scenario, whereas, in the baseline scenario, $23$\% of vehicles have travel time higher than $40$ s. For the scenario shown in Fig. \ref{fig:histTravelTime}, the average travel time has been reduced by $17.62$\%.    

\begin{table}[ht]
\caption{Average travel time of vehicles in the first scenario for the optimal and baseline cases under different traffic volumes.}
\vspace{0.5em}
\centering
\begin{tabular}{c|c|c|c|c} \label{tbl:average travel time}
    Flow & Average number& \multicolumn{2}{c}{Average travel time (s)}& 
    \\
    (veh/h)   & of vehicles &  Baseline              & Optimal          & $\%$
    \\
    \toprule
        $600$  &$44$&  $25.51$  &   $19.41$  & $24$\\
        $800$  &$61$&  $25.14$  &   $20.23$  & $20$\\
        $\num{1000}$  &$76$&  $26.03$  &   $20.59$  & $21$\\
        $\num{1200}$ &$91$&  $26.26$   &   $21.99$  & $16$\\
        $\num{1400}$ &$110$&  $27.27$  &   $24.30$  & $11$\\
\end{tabular}
\end{table}
To investigate the computational complexity of our approach, for each traffic flow in each seed, the time that it takes for a CAV to compute the optimal trajectory is recorded, and then averaged across all the CAVs. For each traffic flow, we choose the seed with a maximum mean of computation time to report. The mean and standard deviation of computation times of CAVs in our optimal proposed framework for different traffic volumes are listed in Table \ref{tbl:computation}. It shows that our approach is computationally feasible and does not grow exponentially with increasing the traffic volume. Please note that since our scheme is decentralized, there is not a relation between the traffic flow and computation times. 

\begin{table}[ht]
\caption{The mean and standard deviation of computation times of CAVs in the first scenario.}
\vspace{0.5em}
\centering
\begin{tabular}{c|c|c|c|c|c} \label{tbl:computation}
    Traffic volume & $600$ & $800$&  $\num{1000}$ & $\num{1200}$& $\num{1400}$
    \\
    \toprule
    Mean (ms) & $0.21$ & $0.19$&  $0.18$ & $0.18$& $0.17$ \\
    Standard deviation (ms) & $0.14$ & $0.14$&  $0.13$ & $0.13$& $0.13$
\end{tabular}
\end{table}
 Our next measure of effectiveness is time-delay, which is computed as a difference between the vehicle's travel time, and the time that it would have taken for the vehicle to cruise with the same speed as the one that it entered the control zone. For CAV $i\in\mathcal{N}(t)$, the time-delay is denoted by $t_{i}^{\text{delay}}$ and given by
 \begin{equation}\label{eq:time-delay}
     t_{i}^{\text{delay}} = (t_i^f - t_i^0)- \dfrac{p_i(t_i^f)-p_i(t_i^0)}{v_i(t_i^0)}.
 \end{equation}

\begin{figure}[ht]
\centering
\includegraphics[width=0.95\linewidth]{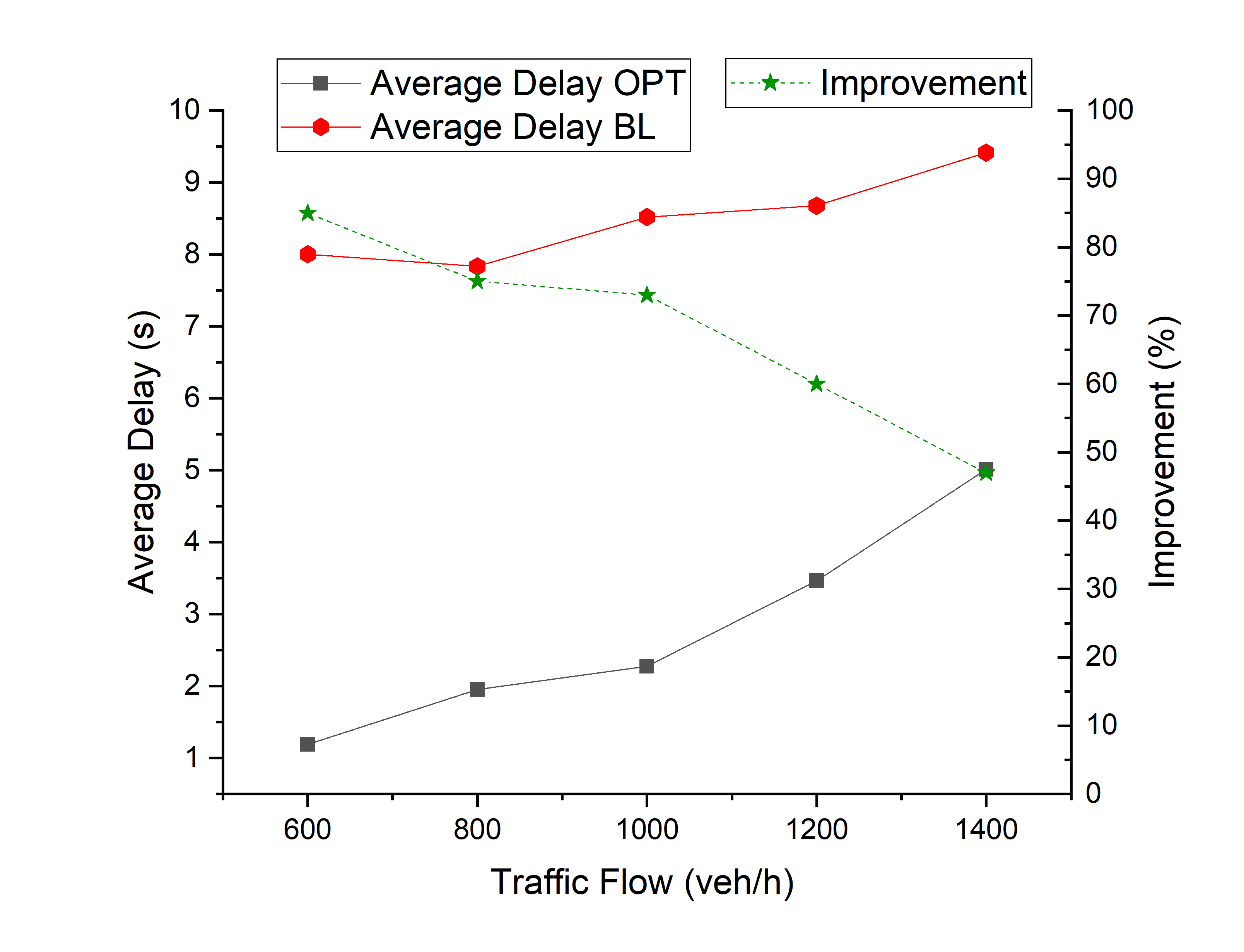}
\caption{Average delay of vehicles in the first scenario} for the baseline and optimal cases.
\label{fig:average delay}%
\end{figure}

 Average delay of all CAVs inside the control zone for the baseline and optimal scenarios at different traffic volumes ranging from $600$ veh/h to $\num{1400}$ veh/h per lane for each entry, along with the percentage of improvement are illustrated in Fig. \ref{fig:average delay}. As shown in Fig. \ref{fig:average delay}, in our proposed framework the average delay has been reduced $47$\% - $85$\% compared to the baseline scenario. 
 
\begin{figure*}[ht]
    \centering
\subfloat[][]{\includegraphics[width=0.33\linewidth]{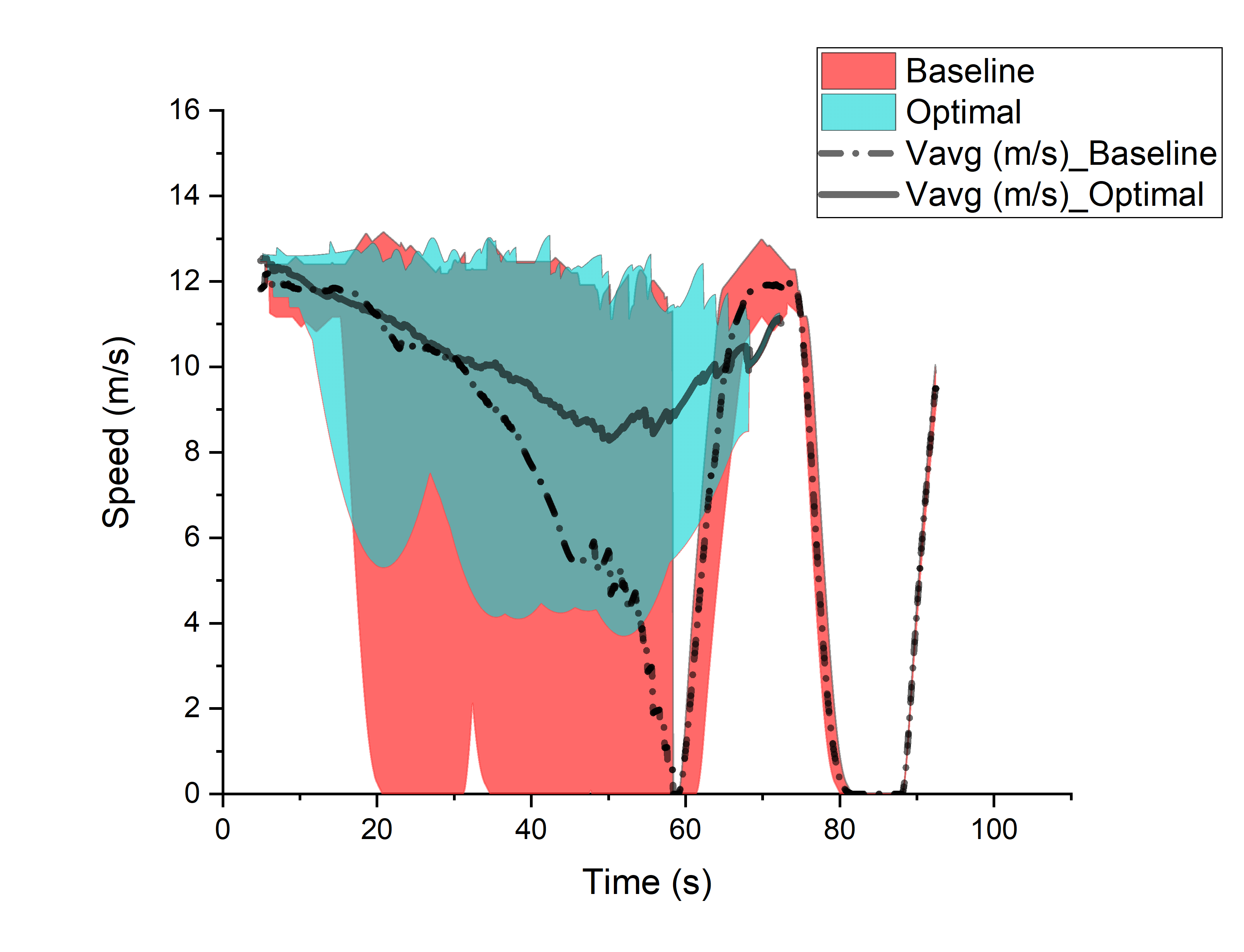}\label{a}}
\subfloat[][]{\includegraphics[width=0.33\linewidth]{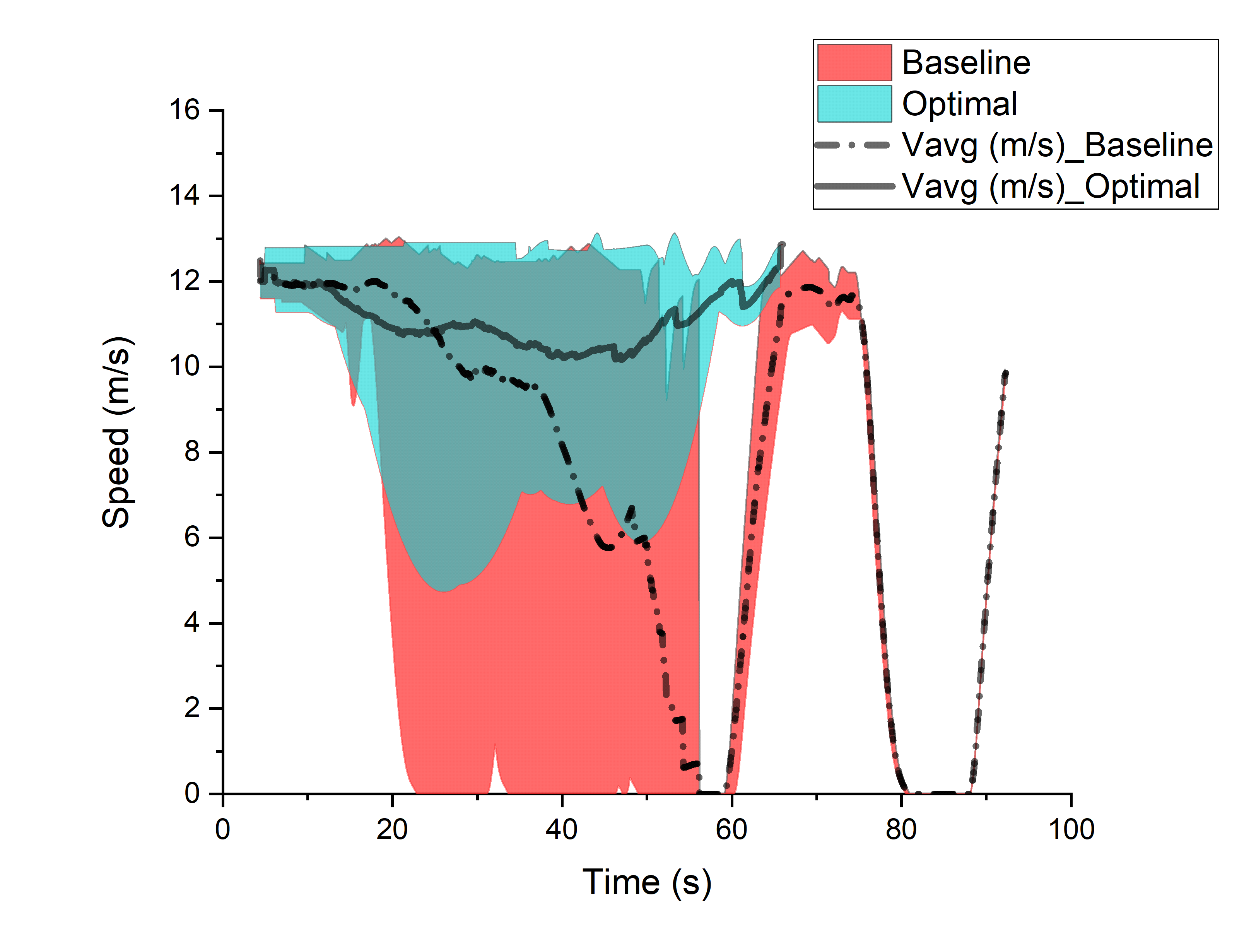}\label{b}}
\subfloat[][]{\includegraphics[width=0.33\linewidth]{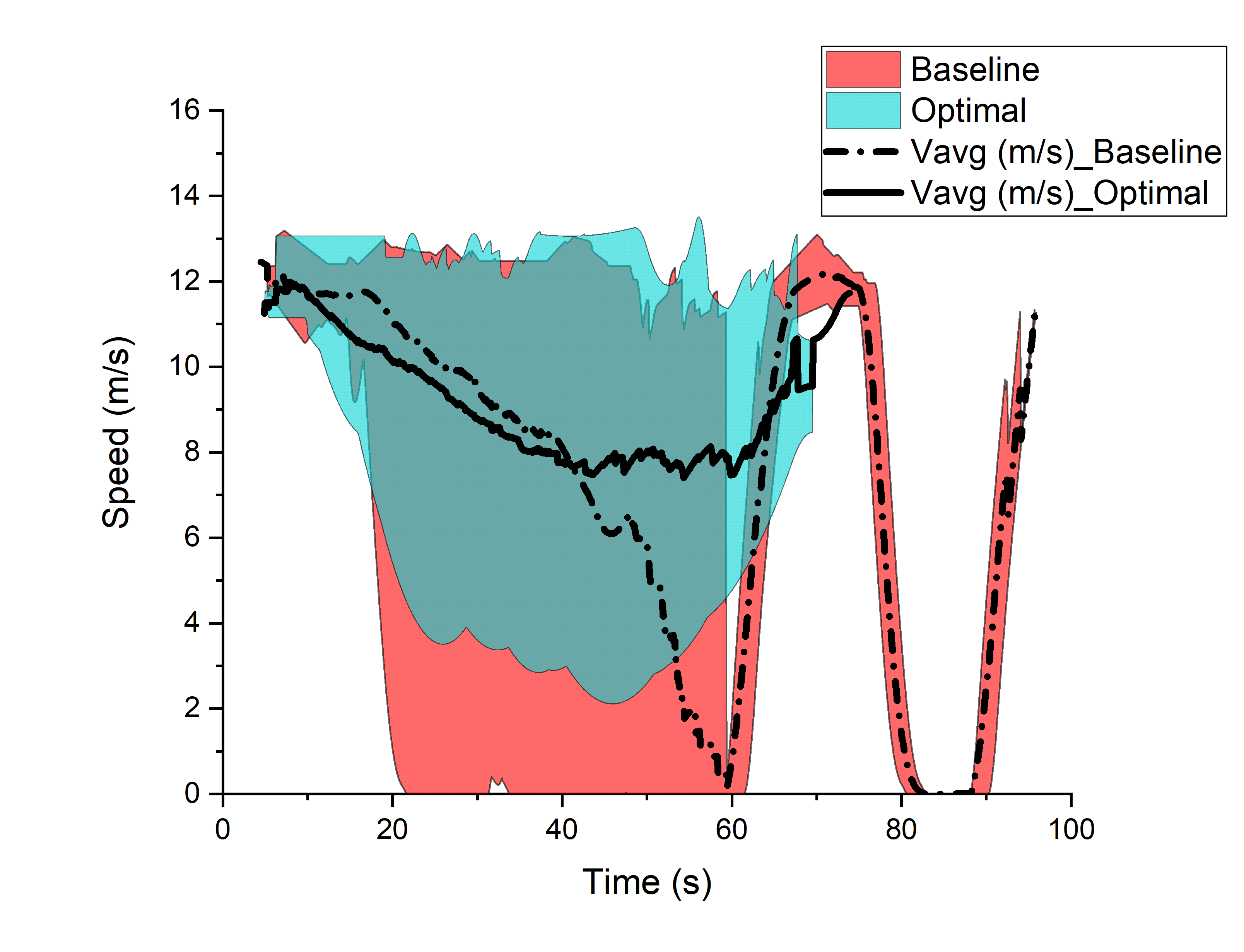}\label{c}}
    \caption{The instantaneous average, maximum and minimum speed of CAVs inside the control zone in the first scenario for the baseline and optimal cases with traffic volume \protect\subref{a} $600$ veh/h, \protect\subref{b} $\num{1000}$ veh/h and \protect\subref{c} $\num{1400}$ veh/h.}
    \label{fig:rangeSpeed}
\end{figure*}

\begin{figure}[ht]
\centering
\includegraphics[width=0.95\linewidth]{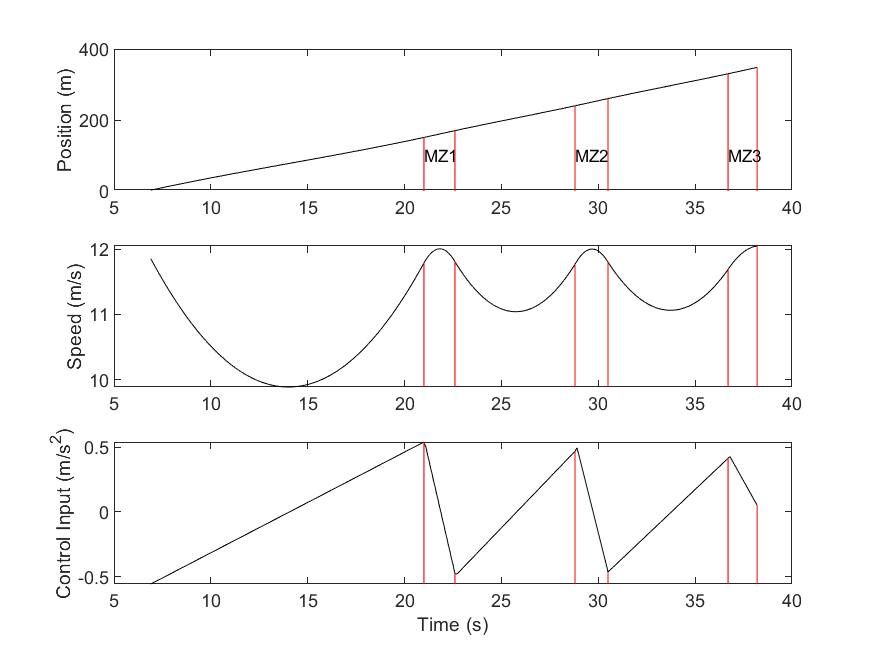}
\caption{The position, speed and control input for a CAV entering the control zone from east for the optimal case in the first scenario.}
\label{fig:profileOptimal}%
\end{figure}

\begin{figure}[ht]
\centering
\includegraphics[width=0.95\linewidth]{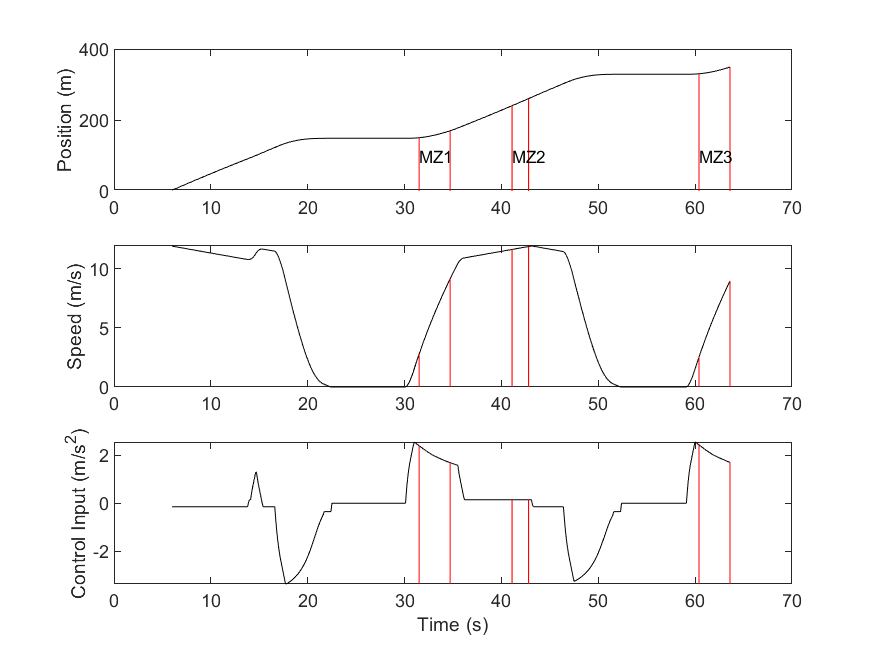}
\caption{The position, speed, and control input for a CAV entering the control zone from east for the baseline case in the first scenario.}
\label{fig:profileBaseline}%
\end{figure}

The instantaneous average, maximum, and minimum speed of CAVs inside the control zone for the baseline and optimal scenarios with traffic volume $600$ veh/h, $\num{1000}$ veh/h, and $\num{1400}$ veh/h for three randomly selected seeds are shown in Fig. \ref{fig:rangeSpeed}. The average speed for the optimal scenario is higher than the average speed in the baseline scenario most of the time, which shows improved traffic throughput. The instantaneous minimum speed for all traffic volumes in the optimal scenario is positive indicating smooth traffic flow, compared to the baseline scenario, which experiences much stopping due to the traffic lights. The position, speed, and control input for a CAV entering the control zone from the east are shown in Fig. \ref{fig:profileOptimal} for the optimal proposed framework, and in Fig. \ref{fig:profileBaseline} for the baseline scenario. As it was shown earlier (Theorem \ref{theorem2}), the control input is a continuous piecewise linear function. We can see in Fig. \ref{fig:profileOptimal} that the CAV needs to accelerate to satisfy the imposed average speed at each merging zone. Choosing the optimal speed for $v_{\text{avg}}^z$ may reduce this oscillation, and potentially improve the system's efficiency.


To evaluate fuel efficiency improvement achieved by our proposed framework, we use a polynomial meta-model proposed in \cite{kamal2012model}, which approximates the fuel consumption in ml/s as a function of speed and control input of a CAV and coefficients obtained from an engine torque-speed-efficiency map of a typical car. Table \ref{tbl:average fuel consumption} summarizes the average fuel consumption and the average cumulative fuel consumption for the optimal and baseline scenarios at different traffic volumes for which five simulations with different random seeds were performed, and the results were averaged. It can be noted that our optimal framework results in better fuel efficiency compared to the baseline scenario. Within our optimal framework, the average cumulative fuel consumption has been improved by $32$\% - $55$\% compared to the baseline scenario with two-phase traffic signals.

\begin{table}[ht]
\caption{Average fuel consumption and average cumulative fuel consumption in first scenario for the optimal and baseline cases under different traffic volumes.}
\vspace{0.5em}
\centering
\begin{tabular}{c|c|c|c|c|c} \label{tbl:average fuel consumption}
    &\multicolumn{2}{c}{Average fuel}&\multicolumn{2}{c}{Average cumulative}&\\
  Flow  & \multicolumn{2}{c}{consumption (ml/s)}& \multicolumn{2}{c}{fuel consumption (ml)}&\\
    (veh/h)   &Baseline & Optimal&Baseline &  Optimal& \%
    \\
    \toprule
        $600$  &$0.32$&  $0.19$  &   $8.87$  & $3.89$&$55$\\
        $800$  &$0.33$&  $0.21$  &   $8.91$  & $4.58$&$48$\\
        $\num{1000}$  &$0.34$&  $0.23$  &   $9.46$  & $4.92$&$48$\\
        $\num{1200}$ &$0.34$&  $0.25$   &  $9.54$  & $5.80$&$39$\\
        $\num{1400}$ &$0.34$&  $0.27$  &   $10.01$  & $6.79$&$32$\\
\end{tabular}
\end{table}

\begin{figure}[ht]
\centering
\includegraphics[width=0.95\linewidth]{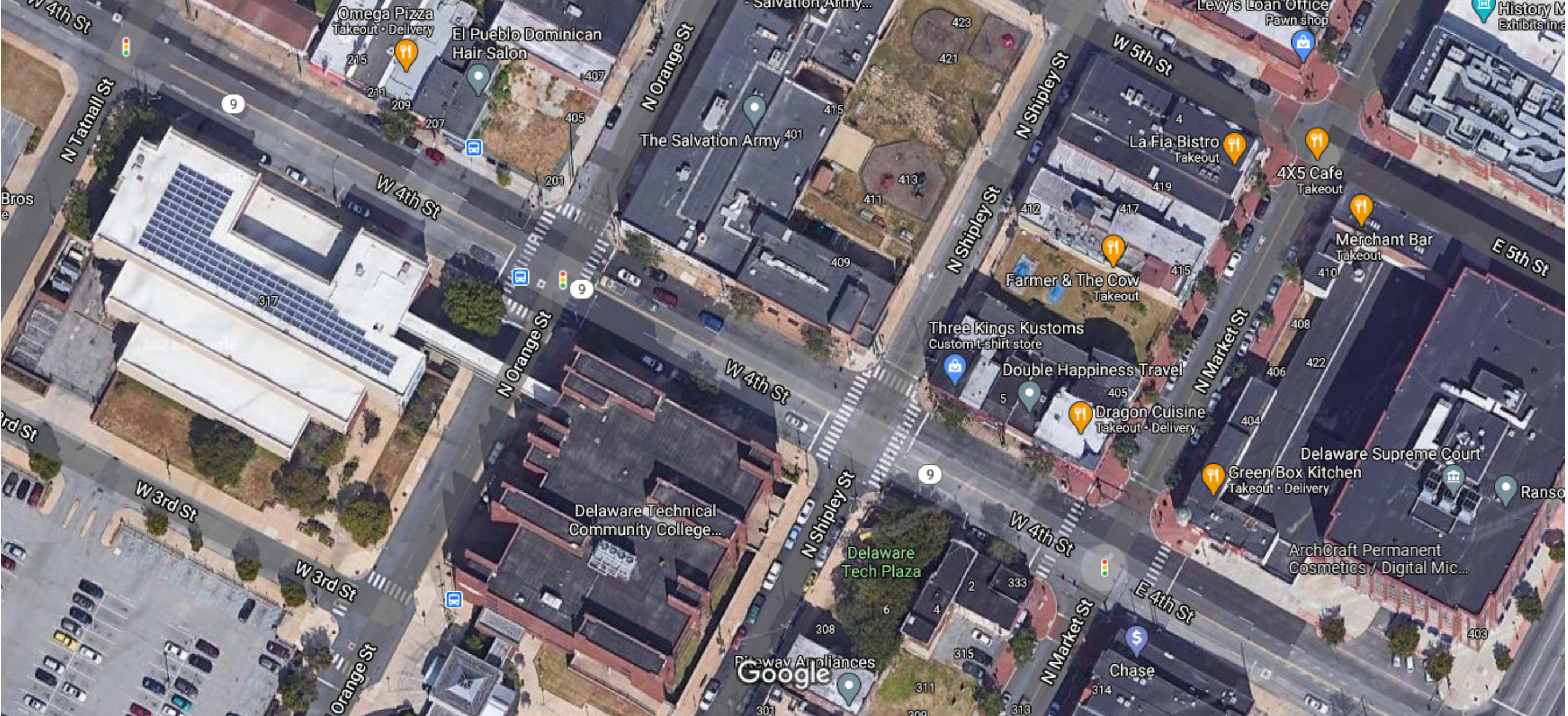}
\caption{An asymmetric corridor in W $4$th street in Wilmington, Delaware, second scenario.}
\label{fig:secondScenario}%
\end{figure}

In the second scenario, we consider an asymmetric corridor in W $4$th street in Wilmington, Delaware consisting of three adjacent intersections with N Orange street, N Shipley street, and N Market street (see Fig. \ref{fig:secondScenario}). We consider the length of each road connecting to the intersections to be $L=150$ m. The vehicles enter the control zone with initial speed uniformly distributed between $8$ m/s to $11$ m/s from each entry with equal traffic volumes. For each traffic volume, we performed five simulations with different random seeds and averaged the results.

Table \ref{tbl:average travel time scenario2} contains the average travel time and average delay of all CAVs inside the control zone for the baseline and optimal cases, respectively, at different traffic volumes ranging from $600$ veh/h to $\num{1400}$ veh/h per lane for each entry. The results indicate that three is $21\% -35\%$ reduction in the average travel time, and $57\% -84\%$ decrease in the average delay compared to the baseline scenario with two-phase traffic signals.

\begin{table}[ht]
\caption{Average travel time and delay of vehicles in in the second scenario for the optimal and baseline cases under different traffic volumes.}
\vspace{0.5em}
\centering
\begin{tabular}{c|c|c|c|c|c|c} \label{tbl:average travel time scenario2}
    Flow & \multicolumn{2}{c}{Avg. travel time (s)}& & \multicolumn{2}{c}{Avg. delay (s)}& 
    \\
    (veh/h)   & Baseline &  Optimal&$\%$& Baseline &  Optimal&$\%$
    \\
    \toprule
        $600$  &  $37.72$  &   $24.53$  & $35$& $15.92$ & $2.56$ & $84$\\
        $800$  & $39.67$  &   $25.89$  & $35$& $17.80$ & $3.47$ & $80$\\
        $\num{1000}$  & $40.10$  &   $28.65$  & $29$& $18.17$ & $5.42$ & $70$\\
        $\num{1200}$ & $40.69$   &   $30.38$  & $25$& $18.57$ & $6.52$ & $65$\\
        $\num{1400}$ &  $42.25$  &   $33.38$  & $21$& $19.98$ & $8.68$ & $57$\\
\end{tabular}
\end{table}

Table \ref{tbl:average fuel consumption second sceario} shows the average fuel consumption and the average cumulative fuel consumption for the optimal and baseline cases at different traffic volumes. The results further support the improvement in fuel efficiency by using the optimal framework. Namely, the average cumulative fuel consumption has been improved by $54$\% - $62$\% compared to the baseline scenario with two-phase traffic signals.
\begin{table}[ht]
\caption{Average fuel consumption and average cumulative fuel consumption in second scenario for the optimal and baseline cases under different traffic volumes.}
\vspace{0.5em}
\centering
\begin{tabular}{c|c|c|c|c|c} \label{tbl:average fuel consumption second sceario}
    &\multicolumn{2}{c}{Average fuel}&\multicolumn{2}{c}{Average cumulative}&\\
  Flow  & \multicolumn{2}{c}{consumption (ml/s)}& \multicolumn{2}{c}{fuel consumption (ml)}&\\
    (veh/h)   &Baseline & Optimal&Baseline &  Optimal& \%
    \\
    \toprule
        $600$  &$0.25$&  $0.14$  &   $9.36$  & $3.86$&$59$\\
        $800$  &$0.25$&  $0.13$  &   $9.74$  & $3.73$&$62$\\
        $\num{1000}$  &$0.24$&  $0.13$  &   $9.79$  & $4.14$&$58$\\
        $\num{1200}$ &$0.24$&  $0.13$   &  $9.96$  & $4.44$&$55$\\
        $\num{1400}$ &$0.24$&  $0.14$  &   $10.18$  & $4.73$&$54$\\
\end{tabular}
\end{table}

\section{Concluding remarks and future research}\label{section7}
In this paper, we proposed a bi-level decentralized coordination framework for CAVs at multiple adjacent multi-lane signal-free intersections closely distanced from each other. In the upper-level planning, each CAV recursively computes the energy-optimal arrival time at each intersection along its path, while ensuring both lateral and rear-end safety. By introducing a lane-changing zone, each CAV investigates the feasibility of a lane-changing maneuver and determines the optimal lane to occupy, aimed at improving the traffic throughput. In the low-level planning, we formulated an optimal control problem for each CAV with the interior-point constraints, the solution of which yields the energy optimal control input (acceleration/deceleration), given the time from the upper-level problem. We developed a recursive structure for the upper-level planning, and also derived an analytical solution for the optimal control problem with interior-point constraints, that can be implemented in real time. In addition, we enhanced our bi-level framework to guarantee safety in the presence of a bounded steady-state error in tracking the positions of CAVs. Finally, our proposed framework exhibited reduction in fuel-consumption, traffic delay, and improvement in the travel time compared to the baseline scenario in different traffic volumes ranging from $600$ veh/h to $\num{1400}$ veh/h for both symmetric and asymmetric adjacent intersections.

Owing to the fact that solving a constrained solution leads to solving a system of non-linear equations that might be hard to solve in real-time for some cases, a different approach has been explored in \cite{Malikopoulos2020} and \cite{Malikopoulos2019CDC}, in which the upper-level optimization problem yields a final time that results in the unconstrained energy-optimal solution in the low-level problem. This approach has also been experimentally validated in \cite{chalaki2020experimental} at University of Delaware's Scaled Smart City for multi-lane roundabouts. 

Ongoing research considers uncertainty in the framework originated from the vehicle-level control \cite{chalaki2021RobustGP} and also investigates the effects of errors and delays in the V2V and V2I communication. Recently, several studies have developed eco-driving approaches for signalized intersections under mixed-traffic scenarios \cite{wang2019cooperative,yang2021cooperative,yang2020eco}. However, coordination for mixed-traffic scenarios and the interaction of human-driven vehicles and CAVs is still an open research question and a potential direction for future research.

%%%%Bibiliography%%%%
\bibliographystyle{IEEEtran} 
\bibliography{bib/IDS_Publications_04222021.bib, bib/ref.bib}

\begin{thebibliography}{10}
\providecommand{\url}[1]{#1}
\csname url@samestyle\endcsname
\providecommand{\newblock}{\relax}
\providecommand{\bibinfo}[2]{#2}
\providecommand{\BIBentrySTDinterwordspacing}{\spaceskip=0pt\relax}
\providecommand{\BIBentryALTinterwordstretchfactor}{4}
\providecommand{\BIBentryALTinterwordspacing}{\spaceskip=\fontdimen2\font plus
\BIBentryALTinterwordstretchfactor\fontdimen3\font minus
  \fontdimen4\font\relax}
\providecommand{\BIBforeignlanguage}[2]{{%
\expandafter\ifx\csname l@#1\endcsname\relax
\typeout{** WARNING: IEEEtran.bst: No hyphenation pattern has been}%
\typeout{** loaded for the language `#1'. Using the pattern for}%
\typeout{** the default language instead.}%
\else
\language=\csname l@#1\endcsname
\fi
#2}}
\providecommand{\BIBdecl}{\relax}
\BIBdecl

\bibitem{united20182018}
{United Nations}, ``2018 revision of world urbanization prospects,'' 2018.

\bibitem{Schrank2019}
B.~Schrank, B.~Eisele, and T.~Lomax, ``{2019 Urban Mobility Scorecard},'' Texas
  A\& M Transportation Institute, Tech. Rep., 2019.

\bibitem{USDOT2}
{National Center for Statistics and Analysis}, ``Police-reported motor vehicle
  traffic crashes in 2018 (traffic safety facts research note),'' Tech. Rep.
  DOT HS 812 860, 2019.

\bibitem{Klein2016a}
I.~Klein and E.~Ben-Elia, ``{Emergence of cooperation in congested road
  networks using ICT and future and emerging technologies: A game-based
  review},'' \emph{Transportation Research Part C: Emerging Technologies},
  vol.~72, pp. 10--28, 2016.

\bibitem{Melo2017a}
S.~Melo, J.~Macedo, and P.~Baptista, ``{Guiding cities to pursue a smart
  mobility paradigm: An example from vehicle routing guidance and its traffic
  and operational effects},'' \emph{Research in Transportation Economics},
  vol.~65, pp. 24--33, 2017.

\bibitem{Levine1966}
W.~Levine and M.~Athans, ``{On the optimal error regulation of a string of
  moving vehicles},'' \emph{IEEE Transactions on Automatic Control}, vol.~11,
  no.~3, pp. 355--361, 1966.

\bibitem{Athans1969}
M.~Athans, ``{A unified approach to the vehicle-merging problem},''
  \emph{Transportation Research}, vol.~3, no.~1, pp. 123--133, 1969.

\bibitem{Dresner2008}
K.~Dresner and P.~Stone, ``{A Multiagent Approach to Autonomous Intersection
  Management},'' \emph{Journal of Artificial Intelligence Research}, vol.~31,
  pp. 591--653, 2008.

\bibitem{Lee2012a}
J.~Lee and B.~Park, ``Development and evaluation of a cooperative vehicle
  intersection control algorithm under the connected vehicles environment,''
  \emph{IEEE Transactions on Intelligent Transportation Systems}, vol.~13,
  no.~1, pp. 81--90, 2012.

\bibitem{gregoire2014priority}
J.~Gregoire, S.~Bonnabel, and A.~De~La~Fortelle, ``Priority-based intersection
  management with kinodynamic constraints,'' in \emph{2014 European Control
  Conference (ECC)}.\hskip 1em plus 0.5em minus 0.4em\relax IEEE, 2014, pp.
  2902--2907.

\bibitem{fayazi2018mixed}
S.~A. Fayazi and A.~Vahidi, ``Mixed-integer linear programming for optimal
  scheduling of autonomous vehicle intersection crossing,'' \emph{IEEE
  Transactions on Intelligent Vehicles}, vol.~3, no.~3, pp. 287--299, 2018.

\bibitem{bichiou2018developing}
Y.~Bichiou and H.~A. Rakha, ``Developing an optimal intersection control system
  for automated connected vehicles,'' \emph{IEEE Transactions on Intelligent
  Transportation Systems}, vol.~20, no.~5, pp. 1908--1916, 2018.

\bibitem{borek2019economic}
J.~Borek, B.~Groelke, C.~Earnhardt, and C.~Vermillion, ``Economic optimal
  control for minimizing fuel consumption of heavy-duty trucks in a highway
  environment,'' \emph{IEEE Transactions on Control Systems Technology}, 2019.

\bibitem{du2018hierarchical}
Z.~Du, B.~HomChaudhuri, and P.~Pisu, ``Hierarchical distributed coordination
  strategy of connected and automated vehicles at multiple intersections,''
  \emph{Journal of Intelligent Transportation Systems}, vol.~22, no.~2, pp.
  144--158, 2018.

\bibitem{Colombo2015}
A.~Colombo and D.~{Del Vecchio}, ``{Least restrictive supervisors for
  intersection collision avoidance: A scheduling approach},'' \emph{IEEE
  Transactions on Automatic Control}, 2015.

\bibitem{Colombo2014}
A.~Colombo, ``{A mathematical framework for cooperative collision avoidance of
  human-driven vehicles at intersections},'' in \emph{2014 11th International
  Symposium on Wireless Communications Systems, ISWCS 2014 - Proceedings},
  2014.

\bibitem{Ahn2014}
H.~Ahn, A.~Colombo, and D.~{Del Vecchio}, ``{Supervisory control for
  intersection collision avoidance in the presence of uncontrolled vehicles},''
  in \emph{Proceedings of the American Control Conference}, 2014.

\bibitem{Ahn2016a}
H.~Ahn and D.~{Del Vecchio}, ``{Semi-autonomous intersection collision
  avoidance through job-shop scheduling},'' \emph{Proceedings of the 19th
  International Conference on Hybrid Systems: Computation and Control - HSCC
  '16}, pp. 185--194, 2016.

\bibitem{DeCampos2015a}
G.~R. {De Campos}, F.~{Della Rossa}, and A.~Colombo, ``{Optimal and least
  restrictive supervisory control: Safety verification methods for human-driven
  vehicles at traffic intersections},'' \emph{Proceedings of the IEEE
  Conference on Decision and Control}, pp. 1707--1712, 2015.

\bibitem{Makarem2012}
L.~Makarem and D.~Gillet, ``{Fluent coordination of autonomous vehicles at
  intersections},'' \emph{2012 IEEE International Conference on Systems, Man,
  and Cybernetics (SMC)}, pp. 2557--2562, Oct. 2012.

\bibitem{wu2014distributed}
W.~Wu, J.~Zhang, A.~Luo, and J.~Cao, ``Distributed mutual exclusion algorithms
  for intersection traffic control,'' \emph{IEEE Transactions on Parallel and
  Distributed Systems}, vol.~26, no.~1, pp. 65--74, 2014.

\bibitem{azimi2014stip}
R.~Azimi, G.~Bhatia, R.~R. Rajkumar, and P.~Mudalige, ``{STIP}: Spatio-temporal
  intersection protocols for autonomous vehicles,'' in \emph{2014 ACM/IEEE
  International Conference on Cyber-Physical Systems (ICCPS)}.\hskip 1em plus
  0.5em minus 0.4em\relax IEEE, 2014, pp. 1--12.

\bibitem{hult2018optimal}
R.~Hult, M.~Zanon, S.~Gros, and P.~Falcone, ``Optimal coordination of automated
  vehicles at intersections: Theory and experiments,'' \emph{IEEE Transactions
  on Control Systems Technology}, vol.~27, no.~6, pp. 2510--2525, 2018.

\bibitem{kim2014mpc}
K.-D. Kim and P.~R. Kumar, ``An {MPC}-based approach to provable system-wide
  safety and liveness of autonomous ground traffic,'' \emph{IEEE Transactions
  on Automatic Control}, vol.~59, no.~12, pp. 3341--3356, 2014.

\bibitem{campos2014cooperative}
G.~R. Campos, P.~Falcone, H.~Wymeersch, R.~Hult, and J.~Sj{\"o}berg,
  ``Cooperative receding horizon conflict resolution at traffic
  intersections,'' in \emph{53rd IEEE Conference on Decision and
  Control}.\hskip 1em plus 0.5em minus 0.4em\relax IEEE, 2014, pp. 2932--2937.

\bibitem{kloock2019distributed}
M.~Kloock, P.~Scheffe, S.~Marquardt, J.~Maczijewski, B.~Alrifaee, and
  S.~Kowalewski, ``Distributed model predictive intersection control of
  multiple vehicles,'' in \emph{2019 IEEE Intelligent Transportation Systems
  Conference (ITSC)}.\hskip 1em plus 0.5em minus 0.4em\relax IEEE, 2019, pp.
  1735--1740.

\bibitem{Malikopoulos2017}
A.~A. Malikopoulos, C.~G. Cassandras, and Y.~Zhang, ``A decentralized
  energy-optimal control framework for connected automated vehicles at
  signal-free intersections,'' \emph{Automatica}, vol.~93, pp. 244--256, 2018.

\bibitem{malikopoulos2019ACC}
A.~A. Malikopoulos and L.~Zhao, ``A closed-form analytical solution for optimal
  coordination of connected and automated vehicles,'' in \emph{2019 American
  Control Conference (ACC)}.\hskip 1em plus 0.5em minus 0.4em\relax IEEE, 2019,
  pp. 3599--3604.

\bibitem{Mahbub2019ACC}
A.~M.~I. Mahbub, L.~Zhao, D.~Assanis, and A.~A. Malikopoulos, ``{Energy-Optimal
  Coordination of Connected and Automated Vehicles at Multiple
  Intersections},'' in \emph{Proceedings of 2019 American Control Conference},
  2019, pp. 2664--2669.

\bibitem{zhang2019decentralized}
Y.~Zhang and C.~G. Cassandras, ``Decentralized optimal control of connected
  automated vehicles at signal-free intersections including comfort-constrained
  turns and safety guarantees,'' \emph{Automatica}, vol. 109, p. 108563, 2019.

\bibitem{zhang2019joint}
Y.~{Zhang} and C.~G. {Cassandras}, ``Joint time and energy-optimal control of
  connected automated vehicles at signal-free intersections with
  speed-dependent safety guarantees,'' in \emph{2019 IEEE 58th Conference on
  Decision and Control (CDC)}, 2019, pp. 329--334.

\bibitem{guanetti2018control}
J.~Guanetti, Y.~Kim, and F.~Borrelli, ``Control of connected and automated
  vehicles: State of the art and future challenges,'' \emph{Annual Reviews in
  Control}, vol.~45, pp. 18--40, 2018.

\bibitem{Rios-Torres2017}
J.~Rios-Torres and A.~A. Malikopoulos, ``A survey on the coordination of
  connected and automated vehicles at intersections and merging at highway
  on-ramps,'' \emph{IEEE Transactions on Intelligent Transportation Systems},
  vol.~18, no.~5, pp. 1066--1077, 2016.

\bibitem{mahbub2020decentralized}
A.~I. Mahbub, A.~A. Malikopoulos, and L.~Zhao, ``Decentralized optimal
  coordination of connected and automated vehicles for multiple traffic
  scenarios,'' \emph{Automatica}, vol. 117, p. 108958, 2020.

\bibitem{chalaki2020TITS}
B.~Chalaki and A.~A. Malikopoulos, ``Time-optimal coordination for connected
  and automated vehicles at adjacent intersections,'' \emph{arXiv preprint
  arXiv:1911.04082}, 2020.

\bibitem{xu2019cooperative}
H.~Xu, Y.~Zhang, L.~Li, and W.~Li, ``Cooperative driving at unsignalized
  intersections using tree search,'' \emph{IEEE Transactions on Intelligent
  Transportation Systems}, 2019.

\bibitem{kamal2014vehicle}
M.~A.~S. Kamal, J.-i. Imura, T.~Hayakawa, A.~Ohata, and K.~Aihara, ``A
  vehicle-intersection coordination scheme for smooth flows of traffic without
  using traffic lights,'' \emph{IEEE Transactions on Intelligent Transportation
  Systems}, vol.~16, no.~3, pp. 1136--1147, 2014.

\bibitem{guney2020scheduling}
M.~A. Guney and I.~A. Raptis, ``Scheduling-based optimization for motion
  coordination of autonomous vehicles at multilane intersections,''
  \emph{Journal of Robotics}, vol. 2020, 2020.

\bibitem{Gregoire2014a}
J.~Gregoire, S.~Bonnabel, and A.~{De La Fortelle}, ``{Priority-based
  intersection management with kinodynamic constraints},'' \emph{2014 European
  Control Conference, ECC 2014}, pp. 2902--2907, 2014.

\bibitem{Malikopoulos2020}
A.~A. Malikopoulos, L.~E. Beaver, and I.~V. Chremos, ``Optimal time trajectory
  and coordination for connected and automated vehicles,'' \emph{Automatica},
  vol. 125, p. 109469, 2021.

\bibitem{Malikopoulos:2013aa}
A.~A. Malikopoulos, ``Stochastic optimal control for series hybrid electric
  vehicles,'' \emph{2013 American Control Conference}, pp. 1189--1194, 2013.

\bibitem{bryson1975applied}
A.~E. Bryson and Y.~C. Ho, \emph{Applied optimal control: optimization,
  estimation and control}.\hskip 1em plus 0.5em minus 0.4em\relax CRC Press,
  1975.

\bibitem{PTV}
\emph{PTV VISSIM 11 User manual}, Planung Transport Verkehr AG,, Karlsruhe,
  Germany, 2018.

\bibitem{Wiedemann1974}
R.~Wiedemann, ``Simulation des strassenverkehrsflusses,'' Ph.D. dissertation,
  Universit{\"a}t Karlsruhe, 1974.

\bibitem{kamal2012model}
M.~A.~S. Kamal, M.~Mukai, J.~Murata, and T.~Kawabe, ``Model predictive control
  of vehicles on urban roads for improved fuel economy,'' \emph{IEEE
  Transactions on control systems technology}, vol.~21, no.~3, pp. 831--841,
  2012.

\bibitem{Malikopoulos2019CDC}
A.~A. Malikopoulos and L.~Zhao, ``Optimal path planning for connected and
  automated vehicles at urban intersections,'' in \emph{Proceedings of the 58th
  IEEE Conference on Decision and Control, 2019}.\hskip 1em plus 0.5em minus
  0.4em\relax IEEE, 2019, pp. 1261--1266.

\bibitem{chalaki2020experimental}
B.~Chalaki, L.~E. Beaver, and A.~A. Malikopoulos, ``Experimental validation of
  a real-time optimal controller for coordination of cavs in a multi-lane
  roundabout,'' in \emph{31st IEEE Intelligent Vehicles Symposium (IV)}, 2020,
  pp. 504--509.

\bibitem{chalaki2021RobustGP}
B.~Chalaki and A.~A. Malikopoulos, ``Robust learning-based trajectory planning
  for emerging mobility systems,'' \emph{arXiv preprint arXiv:2103.03313},
  2021.

\bibitem{wang2019cooperative}
Z.~Wang, G.~Wu, and M.~J. Barth, ``Cooperative eco-driving at signalized
  intersections in a partially connected and automated vehicle environment,''
  \emph{IEEE Transactions on Intelligent Transportation Systems}, vol.~21,
  no.~5, pp. 2029--2038, 2019.

\bibitem{yang2021cooperative}
Z.~Yang, Y.~Feng, and H.~X. Liu, ``A cooperative driving framework for urban
  arterials in mixed traffic conditions,'' \emph{Transportation Research Part
  C: Emerging Technologies}, vol. 124, p. 102918, 2021.

\bibitem{yang2020eco}
H.~{Yang}, F.~{Almutairi}, and H.~{Rakha}, ``Eco-driving at signalized
  intersections: A multiple signal optimization approach,'' \emph{IEEE
  Transactions on Intelligent Transportation Systems}, pp. 1--13, 2020.

\end{thebibliography}
%%%%Bio%%%%

\begin{IEEEbiography}
	[{\includegraphics[width=1.1in,height=1.25in,clip,keepaspectratio]{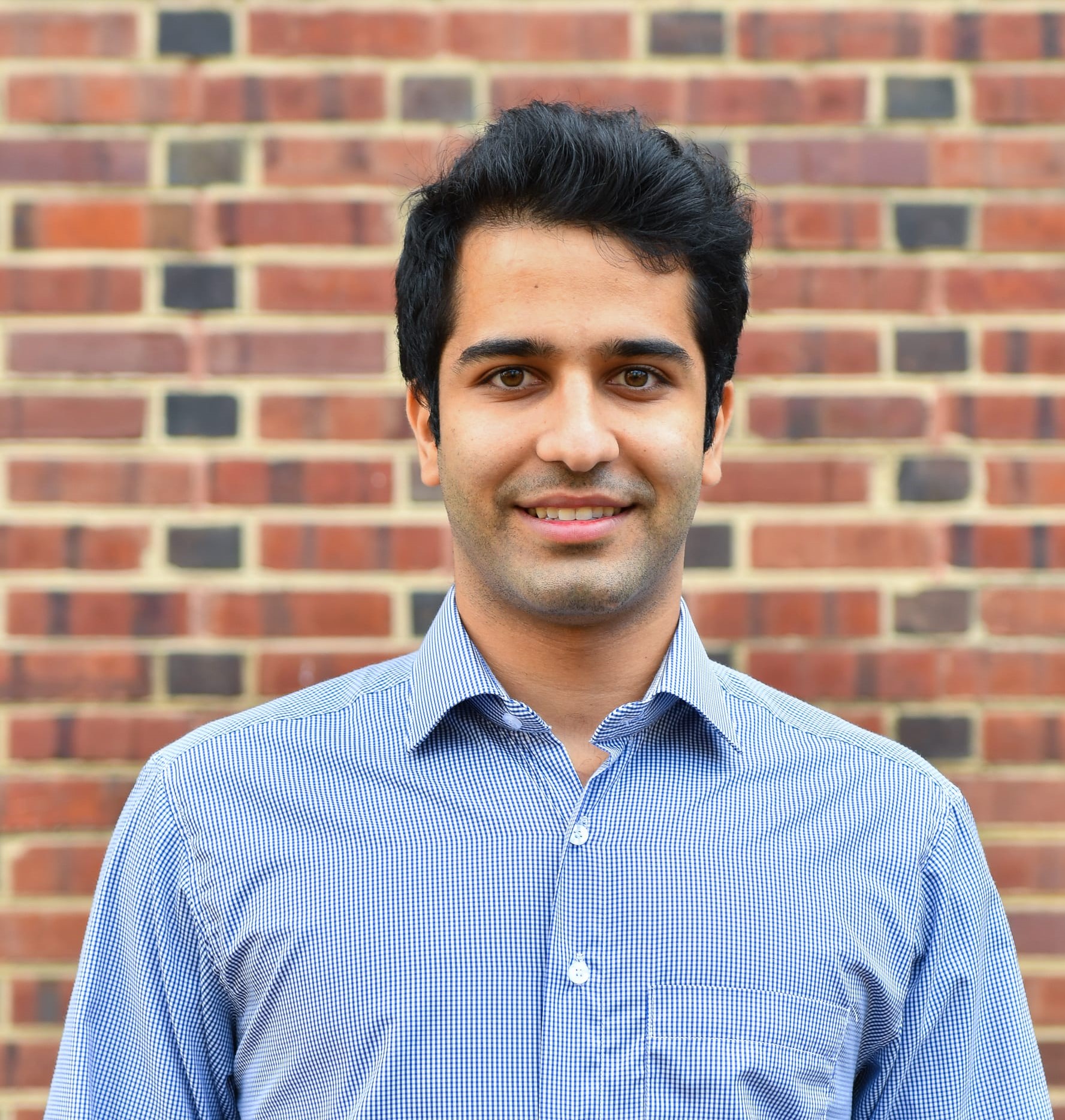}}]
	{Behdad Chalaki} (S’17) received the B.S. degree in mechanical engineering from the University of Tehran, Iran in 2017. He is currently a Ph.D. student in the Information and Decision Science Laboratory in the Department of mechanical engineering at the University of Delaware. His primary research interests lie at the intersections of decentralized optimal control, statistics, and machine learning, with an emphasis on transportation networks. In particular, he is motivated by problems related to improving traffic efficiency and safety in smart cities using optimization techniques. He is a student member of IEEE, SIAM, ASME, and AAAS.
\end{IEEEbiography}
\begin{IEEEbiography}[{\includegraphics[width=1.1in,height=1.25in,clip,keepaspectratio]{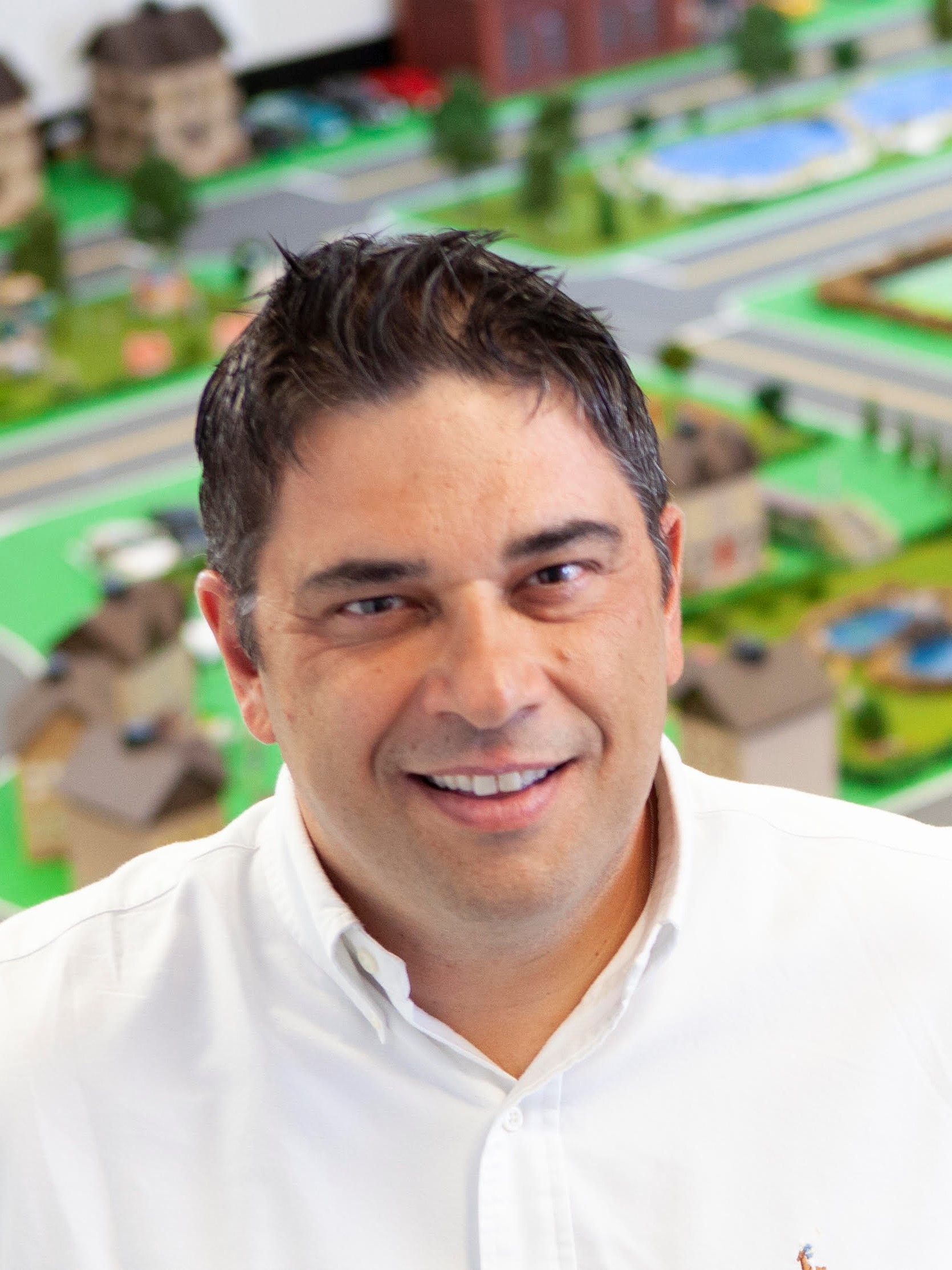}}]{Andreas A. Malikopoulos}
	(S'06--M'09--SM'17) received the Diploma in mechanical engineering from the National Technical University of Athens, Greece, in 2000. He received M.S. and Ph.D. degrees from the department of mechanical engineering at the University of Michigan, Ann Arbor, Michigan, USA, in 2004 and 2008, respectively. 
	He is the Terri Connor Kelly and John Kelly Career Development Associate Professor in the Department of Mechanical Engineering at the University of Delaware (UD), the Director of the Information and Decision Science (IDS) Laboratory, and the Director of the Sociotechnical Systems Center. Before he joined UD, he was the Deputy Director and the Lead of the Sustainable Mobility Theme of the Urban Dynamics Institute at Oak Ridge National Laboratory, and a Senior Researcher with General Motors Global Research \& Development. His research spans several fields, including analysis, optimization, and control of cyber-physical systems; decentralized systems; stochastic scheduling and resource allocation problems; and learning in complex systems. The emphasis is on applications related to smart cities, emerging mobility systems, and sociotechnical systems. He has been an Associate Editor of the IEEE Transactions on Intelligent Vehicles and IEEE Transactions on Intelligent Transportation Systems from 2017 through 2020. He is currently an Associate Editor of Automatica and IEEE Transactions on Automatic Control. He is a member of SIAM, AAAS, and a Fellow of the ASME.
\end{IEEEbiography}

\end{document}